\documentclass[12pt]{article}

\usepackage{booktabs}

\usepackage{algorithm}
\usepackage{algpseudocode}
\usepackage{graphicx,amsmath,amssymb,wrapfig,color,amsfonts}
\usepackage[section]{placeins}
\usepackage{enumitem}
\usepackage{subfig}

\newcommand{\commentout}[1]{}
\newcommand{\keywords}[1]{\textbf{Keywords: } #1}
\DeclareMathAlphabet\mathbfcal{OMS}{cmsy}{b}{n}

\begin{document}
\date{}

\title{A state redistribution algorithm for finite volume schemes
on cut cell meshes}
\author{Marsha Berger\footnote{Courant Institute, New York University, 251 Mercer St.,
NY, NY 10012}  \hspace{1in} Andrew Giuliani $^*$}

\maketitle

\begin{abstract}
In this paper we develop a new technique, called \textit{state redistribution}, 
that allows the use of explicit time stepping when approximating
solutions to hyperbolic conservation laws on embedded boundary grids.  
State redistribution is a postprocessing technique applied after each time 
step or stage of the base finite volume scheme, using a time step
that is proportional to the volume of the full cells.  The idea is to 
stabilize the cut cells by temporarily merging them into larger, 
possibly overlapping neighborhoods, 
then replacing the cut cell values with a stabilized value
that maintains conservation and accuracy.
We present examples of state redistribution using two base schemes:
MUSCL and a second order Method of Lines finite volume scheme. 
State redistribution is used to compute solutions to several 
standard test problems in gas dynamics on cut cell meshes, with both
smooth and discontinuous solutions. We show that our method does not 
reduce the accuracy of the base scheme and 
that it successfully captures shocks in a non-oscillatory manner.
\end{abstract}


\section{Introduction}\label{sec:intro}
Cut cell meshes to solve hyperbolic problems 
are increasingly prevalent due to the ease 
of grid generation for complicated geometries. 
However, these meshes lead to the {\em small cell} problem, which
is still an active area of research, and  
a completely satisfactory solution has not yet been found.
The small cell problem can be explained as follows: explicit
finite volume schemes for hyperbolic problems are subject to a CFL constraint, i.e.,  they typically need to take a time step that is proportional to the mesh width for
stability. However, cut cells can have volumes that are arbitrarily smaller than the regular cells.  This would force the scheme to take overly restrictive time steps, even though the domain is predominantly composed of regular cells that would otherwise determine the stable time step. Special algorithms are needed to prevent this restriction.

The most commonly used stabilization algorithm is called flux
redistribution \cite{chern:colella,vof:colella}. The main idea is illustrated below in two space
dimensions for ease of notation.
The update of the solution average 
$U^n_{i,j}$ on cut cell $(i,j)$  is
\begin{equation*}
V_{i,j} U_{i,j} ^{n+1}  = V_{i,j} U_{i,j}^n  -  \Delta t \sum_{k \in \text{ faces}}A_k \mathbf{F}_k^* \cdot \mathbf{n}_k,
\end{equation*}
Here, $\Delta t$ is the time step,
volume $V_{i,j}$ is the cell volume,
$A_k$, $\mathbf{F}_k^*$, and $\mathbf{n}_k$ are the area, the numerical flux, and the outward facing normal on the $k$th face, respectively.
The above can also be written
\begin{equation*}
V_{i,j} U_{i,j} ^{n+1} = V_{i,j} U_{i,j}^n  +  \delta  M_{i,j} ,
\end{equation*}
where $\delta M_{i,j}$ is the change in the conserved quantity 
on cell $(i,j)$ after one time step.
Instead of using the entire amount of the update in cell $(i,j)$, 
the cut cell only receives a fraction $\eta_{i,j}$ of it.  If the fraction $\eta_{i,j}$
is proportional to the cell's volume fraction $V_{i,j}/V_{\text{full}}$,
where $V_{\text{full}}$ is the volume of an uncut Cartesian cell, the update should be stable. 
To maintain conservation, the rest of the update ($1-\eta_{i,j})\delta M_{i,j}$
is given to the cell's neighbors.  
Flux redistribution has already been implemented for three dimensional
calculations due to its simplicity. However it is only first order accurate at the cut cells.

Cell merging \cite{icm-2003,Chung:2006}  is most frequently the first solution 
that comes to mind for cut cell stabilization. A cut cell is merged with
neighbors until a cell with sufficiently large volume for a stable time
step is obtained. 
It is conceptually simple, but 
we are not aware of any production codes that implement this in a fully
general, robust manner for complicated engineering geometries. 
The $h$-box method \cite{mjb-hel-rjl:hbox2,mjb-hel:hboxsimple}
is a second order accurate method at the cut cells. It extends the 
domain of dependence for the fluxes around a small cell in a 
special way that maintains stability by means of a cancellation
property. It  has not been extended to
three dimensions due to its complexity. 

A newer variation of cell merging is called cell linking \cite{cecereGiacomazzi,
KirkpatrickEtAl:2003, HuKhooAdamsHuang:2006,Chung:2006}.
This has simpler data structures and maintains the original grid. 
In \cite{BalajiMenon:2016}, the authors improve the accuracy of cell linking,
with a third order accurate approach for viscous flow,  and fourth order for 
inviscid flow. 
Their version of cell linking uses a cluster of cells, while still
maintaining each cell in the mesh.  A high order polynomial is fit to
the cluster, and replaces the solution values in the individual cells.
Our state redistribution algorithm has a similar spirit to this, though the 
details are very different. 

Two other approaches in the literature include the use of an implicit
scheme on the cut cells combined with an explicit scheme elsewhere
\cite{May-Berger:JSC}, and a novel flux interpolation scheme which has
the added advantage of being dimensionally split, so easier to
implement \cite{gokhaleNikosKlein:2018}. Our new approach is rather
different from these, however, and we do not pursue these directions
further.

In this paper we propose a stabilization algorithm in
the spirit of flux redistribution. Similar to flux redistribution, state
redistribution is applied as a postprocessing step and is simple to implement,
for the second order accurate case. We perform an unstable update on all 
cells with a fixed $\Delta t$ using a base finite volume scheme, followed
by a postprocessing step based on the conserved state variables, not on the
fluxes.  It is for this reason that we call it state redistribution
(SRD).  Our
approach is fully conservative and can be generalized to high order accuracy,
albeit with  more complexity. The key insight over cell merging was to
recognize that we could maintain both conservation and accuracy using a
weighted convex combination of solution values that takes into account the 
number of overlapping neighborhoods on each cell in the base cut cell
grid.

The important difference between state redistribution and cell mering is
that SRD supports overlapping neighborhoods, and cell merging does not.
As a result, cell merging can be difficult to implement robustly
in three dimensions, since there are many different, possibly
incompatible ways to create non-overlapping merged cells.  State
redistribution does not suffer from this difficulty, and its extension to
three dimensions is straightforward.

The next section illustrates state redistribution  in one space dimension
on a model problem, to give a more intuitive idea without all the details and 
notation of the second order accurate  case.
Section \ref{sec:basefv} discusses the evolution schemes on cut cell meshes
to which SRD is applied.
Section \ref{sec:srdAlg} describes the second order accurate
algorithm in two space dimension.  Computational experiments with the
two-dimensional Euler equations  are presented in
Section \ref{sec:compResults}, and conclusions in Section \ref{sec:conc}. 
We see no reason that this algorithm cannot be extended to higher order
accuracy, and have already started implementing the third and fourth
order accurate cases.

\section{State redistribution in one dimension} \label{sec:srd1d}
We begin this section  by reminding the reader that cell merging can
be written as a postprocessing step. We then show that the extension of
cell merging for overlapping cells does not maintain conservation. 
This motivates our state redistribution approach, which we contrast with
cell merging on a simple one dimensional example. 
Although the small cells do not mimic the cut cells at the 
boundary in higher dimensions, this is still a useful model problem. 

\begin{figure}[h]
\centering
\vspace*{.2in}
\mbox{
\includegraphics[width=0.4\linewidth,trim=30 0 20 0,clip]{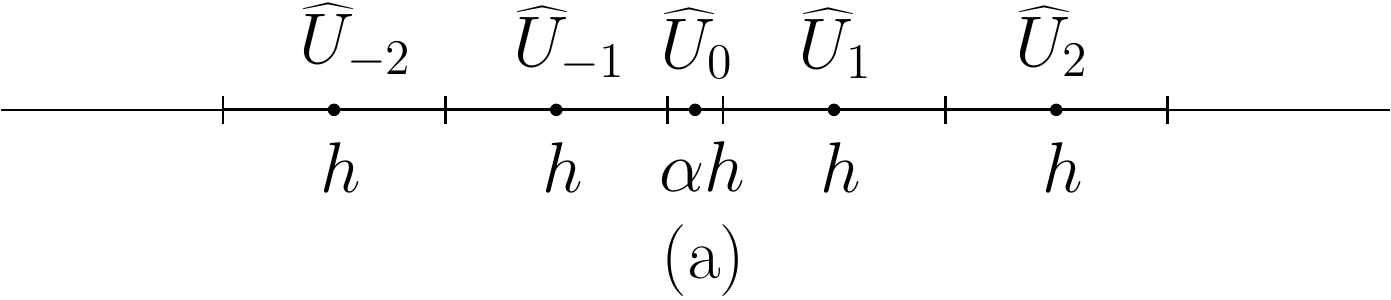} 
\hspace*{.1in}
\includegraphics[width=0.5\linewidth,trim=15 0 25 0,clip]{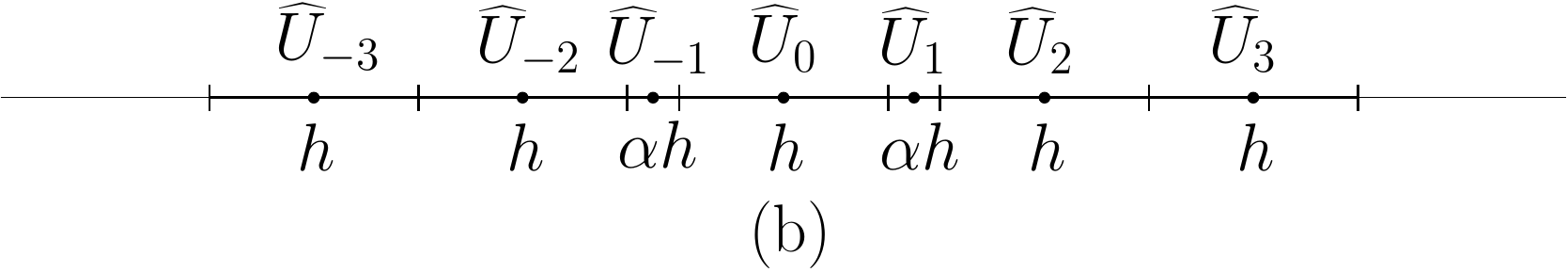}
}
\vspace*{.1in}
\caption{\sf Model problem in one space dimension on nonuniform base grid.
On the left there is one small cell at $i=0$. The right grid has two small 
cells  at $i = -1$ and $1$.  The small and large cell 
sizes are $\alpha h$ ($0<\alpha < 1$), and $h$, respectively. \label{fig:ng1}}
\end{figure}

For the examples in this section we will solve the linear advection equation
\begin{equation}\label{eq:hpde}
u_t + au_x = 0, \quad a>0
\end{equation}
on the nonuniform grid, called the base grid. Equation \eqref{eq:hpde} is
discretized  using the  first order
accurate upwind scheme 
\begin{equation}\label{eq:unstable1d}
\widehat{U}_i = U^n_i - \frac{a \Delta t} {h_i} \, (U^n_i -U^n_{i-1}),
\end{equation}
where $U^n_i$ is the solution average on cell $i$ at time $t^n$, the time
step $\Delta t$ is constant for all cells. On full cells $h_i = h$, and on
the small cells  $h_i = \alpha h$ for $ 0 < \alpha < 1$.
We first use
the grid in Figure \ref{fig:ng1}(a) with one small cell at $i=0$, then the grid
in Figure \ref{fig:ng1}(b) with two small cells at $i = -1$ and $1$.

On the grid in Figure \ref{fig:ng1}(a), cell merging might 
group cells -1,0, and 1 together into a larger, merged cell.
After the unstable update in \eqref{eq:unstable1d}, the  volume-weighted
merged cell average $\widehat{Q}_0$ is computed:
\begin{equation}
\widehat{Q}_0 = \frac{\widehat{U}_{-1} + \alpha \,  \widehat{U}_0 + 
\widehat{U}_1}{2+\alpha} .
\end{equation}
The cells comprising the merged cell
are then replaced by $\widehat{Q}_0$:
\begin{equation}
\widehat{U}_{-1}^{n+1} = \widehat{U}_0^{n+1} = \widehat{U}_1^{n+1} =
\widehat{Q}_0
\end{equation}
This is easily seen to be conservative by checking that
\begin{equation}
\underbrace{h \, (U_{-1}^{n+1} + \alpha \, U_0^{n+1} +
U_1^{n+1})}_{\text{mass after redistribution}}
= (2+\alpha) h \,  \, \widehat{Q}_0 =
\underbrace{ h \, (\widehat{U}_{-1} + \alpha \, \widehat{U}_{0} +
\widehat{U}_{1})}_{\text{mass before redistribution}}
\end{equation}
which equals the values at time $t^n$  except for the mass entering and
leaving this region.

Next, consider the more complicated case of Figure \ref{fig:ng1}(b).
Five cells (indexed by $-3$, $-2$, $0$, $2$, $3$) are large 
with size $h$ and the remaining two cells (indexed by $-1$ and $1$) 
are small with size $\alpha h$.
A first approach to cell merging on this grid
might be to make two merged cell averages, $\widehat{Q}_{-1}$
comprising cells ${-2},{-1},0$, and $\widehat{Q}_1$ comprising cells
$0, 1$ and $2$, associated respectively with the small cells $-1$ and $1$.
The first order version here would assign cells as before, except that
since cell 0 belongs to two neighborhoods, it seems  reasonable
to assign $U_0^{n+1} = \frac{1}{2} (\widehat{Q}_{-1} + \widehat{Q}_1).$  
To check for conservation, we again compute the sum
\begin{equation}
\begin{split}
\underbrace{h \, (U_{-2}^{n+1} + \alpha\,U_{-1}^{n+1} + U_0^{n+1} + \alpha \, U_1^{n+1}
+ U_2^{n+1})}_{\text{mass after redistribution}} = \\[.1in]
h \left [(1+\alpha +1/2) \, \widehat{Q}_{-1} +
 (1/2+\alpha +1) \, \widehat{Q}_{1}\right ]  = \\[.1in]
h \,  (3/2+\alpha) \, \left [ \frac{(\widehat{U}_{-2} + \alpha \widehat{U}_{-1} +
\widehat{U}_0)}{2+\alpha} +
        \frac{(\widehat{U}_{0} + \alpha \widehat{U}_{1} 
        + \widehat{U}_2)}{2+\alpha} \right ]  \neq \\[.08in]
\underbrace{h \,  ( \widehat{U}_{-2} + \alpha \widehat{U}_{-1} +
\widehat{U}_0 +  \alpha \widehat{U}_{1} + 
\widehat{U}_2 )}_{\text{mass before redistribution}}.  
\end{split}
\end{equation}
So at least this extension of cell merging to overlapping cells is not
conservative.

This motivates the state redistribution
procedure, which allows for overlapping cells and stabilizes 
$\widehat{U}_i$ in a conservative manner.
This is done by temporarily merging cells of the grid into 
larger, possibly overlapping, neighborhoods using a specially weighted 
convex combination, and recombining these averages back onto the
grid in a particular fashion.  These merged cells are 
constructed once during a mesh preprocessing step before time stepping.

\begin{figure}[h]
\begin{center}
\vspace*{.1in}
\includegraphics[width=0.75\linewidth]{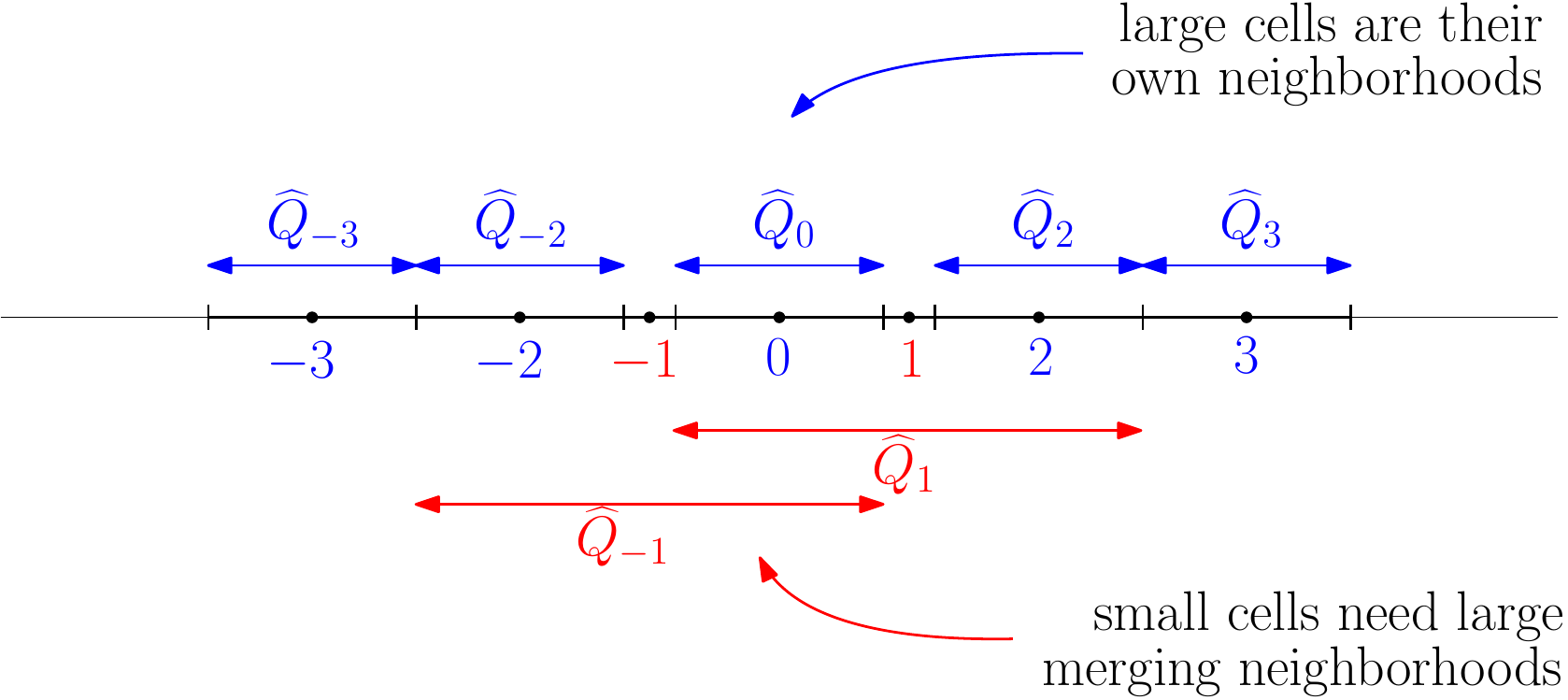} 
\caption{\sf 
The blue arrows indicate the merging neighborhoods associated with the large 
cells $-3$, $-2$, $0$, $2$, and $3$, which are their own neighborhoods.  
The red arrows indicate the neighborhoods of small cells $-1$ and $1$, which
have temporarily merged with their left and right neighbors.  
\label{fig:mn1}}
\end{center}
\end{figure}

\subsubsection*{State redistribution preprocessing}
Each cell in the base grid (both large and small) has a merging neighborhood
associated with it.
This is  a set of neighboring cells with which to temporarily merge.
Neighborhoods share the same index as the cell that generated it.
Small cells merge with their neighbors until the volume of the merging neighborhood is greater than
a threshold, taken here to be half the large cell size of $h/2$.
This is illustrated by the red arrows in Figure \ref{fig:mn1} where the merging neighborhood of 
small cell $-1$ consists of cells $-2$, $-1$, and $0$, and the merging neighborhood of 
small cell 1 consists of cells $0$, $1$, and $2$. 
Note that   both neighborhoods overlap on cell 0. Allowing
for overlaps makes this temporary merging process simpler.
A large cell does not need to merge with neighbors (since $h > h/2$), thus its merging neighborhood is only 
composed of itself.  This is illustrated by the blue arrows in Figure
\ref{fig:mn1}. For example,  the 
merging neighborhood of large cell $-3$ is composed only of itself.

Finally, each cell in the base grid counts the number of neighborhoods that overlap it. 
Cell $-2$ is overlapped by two neighborhoods, indexed by $-1$ and $-2$.
Cell 0 has 3 such neighborhoods, its own, and one from each small cell adjacent to it.

\subsubsection*{State redistribution postprocessing}
Using the above information we can now stabilize \eqref{eq:unstable1d} using the state redistribution 
method on the grid in Figure \ref{fig:ng1}b.  
On each merging neighborhood, we compute a weighted solution average $\widehat Q_i$, where $i$ is the 
index of the merging neighborhood. 
$\widehat Q_i$ is computed by a convex combination of the averages $\widehat{U}_i$ of cells contained in the 
merging neighborhood, weighted by the inverse of its overlap count. 
For example, on merging neighborhood $-1$, the weighted solution average is
\begin{equation}\label{eq:neigh2}
\widehat{Q}_{-1} = \frac{1}{\underbrace{h/2 + \alpha h + h/3}_{\text{weighted volume}}}\biggr( \underbrace{\frac{h}{2} \widehat{U}_{-2} + \alpha h \widehat{U}_{-1} + \frac{h}{3}\widehat{U}_{0}}_{\text{weighted mass}} \biggr).
\end{equation}
In formula \eqref{eq:neigh2}  for the weighted mass,
the cell volume is divided by the number of neighborhoods that overlap the associated cell 
in the base grid. 
For example, the multiplier in front of $\widehat{U}_{-2}$ is $\frac{h}{2}$ since there are 
two neighborhoods (from cells $-1$ and $-2$) that overlap cell $-2$.  
Similarly, the multiplier in front of $\widehat{U}_{-1}$ is $\alpha h$ since there is only 
one neighborhood (its own) that overlaps cell $-1$.
Finally, the multiplier in front of $\widehat{U}_{0}$ is $\frac{h}{3}$ since there are three neighborhoods that overlap cell $0$, i.e., cell $0$ is overlapped by neighborhoods $-1$, $1$, and $0$.
These multipliers are then divided by the weighted volume,  $(h/2 + \alpha h + h/3)$.
i
The weighted solution average on merging neighborhood $1$ is similarly defined as
\begin{equation}\label{eq:neigh3}
\widehat{Q}_{1} = \frac{1}{h/2 + \alpha h + h/3}\left( \frac{h}{2} \widehat{U}_{2} + \alpha h \widehat{U}_{1} + \frac{h}{3}\widehat{U}_{0} \right).
\end{equation}
The weighted solution averages on merging neighborhoods that contain only one cell are simply 
\begin{equation}\label{eq:neigh1}
\widehat{Q}_i = \widehat{U}_i \quad \text{ for } i = -3,-2,0,2,3.
\end{equation}

The stabilized solution average at time $t^{n+1}$ on a cell in the base grid is then 
given by the average of all the weighted neighborhood averages that overlap it.  
On the cell overlapped by three neighborhoods we have
\begin{equation} \label{eq:threeneigh}
U^{n+1}_{0} = \frac{1}{3}(\widehat{Q}_{-1}+\widehat{Q}_{0}+\widehat{Q}_{1}).
\end{equation}
On the  cells overlapped by two neighborhoods, we have
\begin{equation} \label{eq:twoneigh}
U^{n+1}_{-2} = \frac{1}{2}(\widehat{Q}_{-1}+\widehat{Q}_{-2}) \text{ and } U^{n+1}_{2} = \frac{1}{2}(\widehat{Q}_{1}+\widehat{Q}_{2}).
\end{equation}
Finally, on cells overlapped by only one neighborhood,  we have
\begin{equation} \label{eq:oneneigh}
	U^{n+1}_i = \widehat{Q}_i \text{ for } i = -3,-1,1,3.
\end{equation}

We can write the final solution update on the small cells after SRD  at $t^{n+1}$ in terms of the solution 
averages at $t^{n}$, giving
\begin{equation}
\begin{aligned}
U^{n+1}_{-1} &= \frac{2-2\lambda}{5+6\alpha}U^n_0 + \frac{6\alpha - 4 \lambda}{5+6\alpha}U^n_{i-1}+ \frac{3\lambda + 3}{5+6\alpha}U^n_{i-2}+\frac{3\lambda }{5+6\alpha}U^n_{i-3}, \\
U^{n+1}_{1} &= \frac{2+4\lambda}{5+6\alpha}U^n_0 + \frac{6\alpha - 3 \lambda}{5+6\alpha}U^n_{i+1}+ \frac{3-3\lambda}{5+6\alpha}U^n_{i+2}+\frac{2\lambda }{5+6\alpha}U^n_{i-1}.
\end{aligned} \label{eq:finalupdate}
\end{equation}
Before application of the state redistribution method, the weights that multiply the solution averages at time $t^n$ in the base scheme \eqref{eq:unstable1d} become unbounded as $\alpha \rightarrow 0$.
However, after state redistribution this is no longer the case for the weights 
in \eqref{eq:finalupdate}.  This hints at the stability of our modified scheme.  

The state redistribution algorithm allows us to take full time steps as if there were no 
small cells in the grid.  However, it can be seen that the multipliers of $U^n_{i-1}$ and $U^n_{i+1}$ 
in \eqref{eq:finalupdate} are negative when $\alpha$ is small enough.  This means that our scheme is not 
monotone and thus not total variation diminishing.  We note that this is also the case 
with flux redistribution, which has been successfully used in higher dimensions and more
complicated problems.  
The computational examples in Section \ref{sec:compResults} show that
this is not a significant issue for SRD as well.
The advantage of state redistribution  is that it is linearity
preserving, and flux redistribution is not. 

\subsubsection*{Conservation}
We now show that our modified scheme \eqref{eq:threeneigh}, \eqref{eq:twoneigh}, \eqref{eq:oneneigh} conserves mass.
For the portion of the grid in question, the total mass after state redistribution is
\begin{equation}\label{eq:tm1}
	\sum_{i} h_i U^{n+1}_i  = h U^{n+1}_{-3} + h U^{n+1}_{-2} + \alpha h U^{n+1}_{-1} +h U^{n+1}_0+\alpha h U^{n+1}_{1} + h U^{n+1}_{2} + h U^{n+1}_{3},
\end{equation}
where $h_i$ is the local cell size.
Substituting expressions for the final update \eqref{eq:threeneigh}, \eqref{eq:twoneigh}, \eqref{eq:oneneigh} into \eqref{eq:tm1}, we obtain
\begin{equation}\label{eq:tm2}
\begin{aligned}
\sum_{i} h_i U^{n+1}_i  &= h \widehat{Q}_{-3} + \frac{h}{2}(\widehat{Q}_{-1}+\widehat{Q}_{-2}) \\
&+ \alpha h \widehat{Q}_{-1} +h \frac{1}{3}(\widehat{Q}_{-1}+\widehat{Q}_{0}+\widehat{Q}_{1})+\alpha h \widehat{Q}_{1} \\
&+ h \frac{1}{2}(\widehat{Q}_{1}+\widehat{Q}_{2}) + h \widehat{Q}_{3}.
\end{aligned}
\end{equation}
Grouping terms in \eqref{eq:tm2}, we have
\begin{equation}\label{eq:tm3}
\begin{aligned}
\sum_{i} h_i U^{n+1}_i  &= h \widehat{Q}_{-3} + \frac{h}{2}\widehat{Q}_{-2} \\
&+ \left(\frac{h}{2}+\alpha h + \frac{h}{3}\right) \widehat{Q}_{-1} + \frac{h}{3} \widehat{Q}_0 + \left(\frac{h}{2}+\alpha h + \frac{h}{3}\right) \widehat{Q}_{1} \\
&+ \frac{h}{2}\widehat{Q}_{2} + h \widehat{Q}_3.
\end{aligned}
\end{equation}
Substituting the expressions for the neighborhood averages \eqref{eq:neigh3}, \eqref{eq:neigh2}, \eqref{eq:neigh1} into \eqref{eq:tm3}, we obtain
\begin{equation}\label{eq:tm4}
\begin{aligned}
\sum_{i} h_i U^{n+1}_i  &= h \widehat{U}_{-3} + \frac{h}{2}\widehat{U}_{-2} \\
&+ \left(\frac{h}{2} \widehat{U}_{-2} + \alpha h \widehat{U}_{-1} + \frac{h}{3}\widehat{U}_{0}\right) + \frac{h}{3} \widehat{U}_0 + \left(\frac{h}{2} \widehat{U}_{2} + \alpha h \widehat{U}_{1} + \frac{h}{3}\widehat{U}_{0}\right) \\
&+ \frac{h}{2}\widehat{U}_{2} + h \widehat{U}_3.
\end{aligned}
\end{equation}
Simplifying \eqref{eq:tm4}, the mass after state redistribution becomes
\begin{equation}\label{eq:tm5}
\begin{aligned}
\sum_{i} h_i U^{n+1}_i  &= h \widehat{U}_{-3} + h \widehat{U}_{-2} + \alpha h \widehat{U}_{-1} +h \widehat{U}_0+\alpha h \widehat{U}_{1} + h \widehat{U}_{2} + h \widehat{U}_{3},\\
&= \sum_{i} h_i \widehat{U}_i.
\end{aligned}
\end{equation}
Thus, the mass on the grid before and after state redistribution does not change.
Since the base scheme \eqref{eq:unstable1d} is conservative, it follows from
\eqref{eq:tm5} that our modified scheme \eqref{eq:threeneigh},
\eqref{eq:twoneigh}, \eqref{eq:oneneigh} is too.

The stabilized finite volume method \eqref{eq:threeneigh},
\eqref{eq:twoneigh}, \eqref{eq:oneneigh} is first order accurate in space
and time.  In this work, we provide a framework to generalize the state
redistribution method to two dimensional cut cell grids and to 
second order accuracy in space and time.
We will demonstrate with numerical examples that the maximum stable time step is 
not restricted by the small cells, and that the state redistribution method 
is conservative.

\section{Second-Order Accurate Base Schemes}\label{sec:basefv}
\begin{figure}
\begin{center}
\includegraphics[width=3.0in]{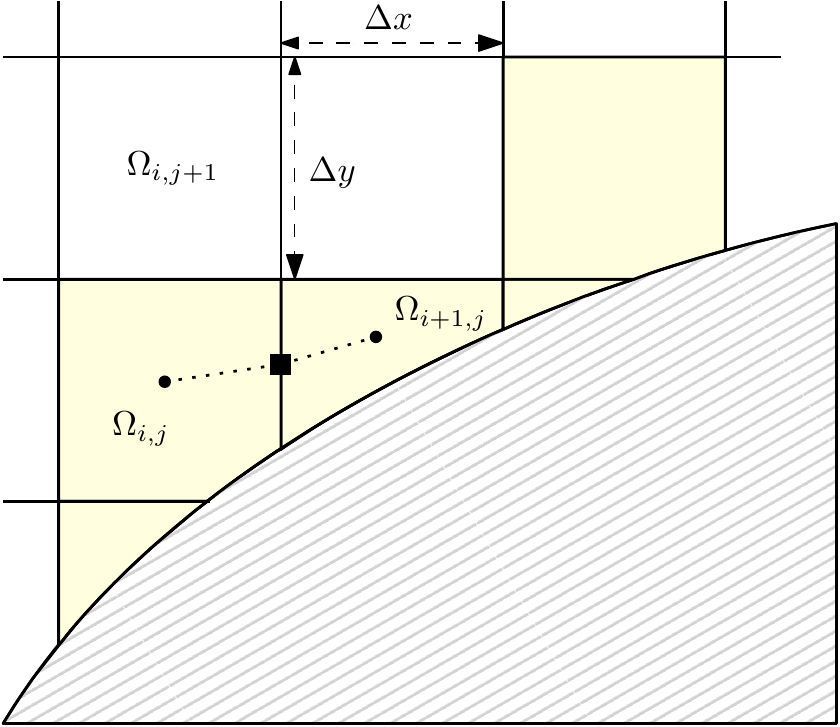}
\caption{\sf Example base grid in two space dimensions. The cells shaded
in yellow are the cut cells.  Both the method of lines (Section
\ref{sec:mol}) and MUSCL (Section \ref{sec:muscl}) schemes require
gradient information to reconstruct to the edge midpoint, indicated with
a square $(\blacksquare)$.} 
\label{fig:2dfig}
\end{center}
\end{figure}

We are interested in solving hyperbolic conservation laws 
\begin{equation}
\begin{aligned} \label{eq:conslaw2D}
\frac{\partial}{\partial t}	\mathbf{u} + \nabla \cdot
\mathbfcal{F}(\mathbf{u})  = \mathbf{0}
\end{aligned}
\end{equation}
on the domain $\Omega \subset \mathbb{R}^2$ where $\mathbf{u}(x,y,t) \in \Omega \times [0,T]$ is a vector of conserved quantities, $T$ is the final time, and $\mathbfcal{F} = [\mathbf{F}, \mathbf{G}]$ is the flux function.  We discretize $\Omega$ into a cut cell mesh of cells $\Omega_{i,j}$.  A typical cut cell mesh, called the base grid, is given in Figure \ref{fig:2dfig}.  On the domain interior, $\Omega_{i,j}$ are Cartesian cells (quadrilaterals) of size $\Delta x$ in the $x$ direction and $\Delta y$ in the $y$ direction.  On the domain boundary there is a border of irregular polygonal cells, called cut cells.  
We use cell-centered discretizations in this work.

There are two issues when applying an explicit finite volume scheme to a cut cell
mesh.  For accuracy, the scheme needs to be modified in the cut cells. Second, 
the scheme needs to be stabilized in the cut cells if using a fixed timestep $\Delta t$ 
based on the full cells. In the $h$-box method these two concerns were addressed
simultaneously, but in general they aren't.

We will use two very different second order discretizations of \eqref{eq:conslaw2D}: the method of lines
(MOL) approach and the MUSCL scheme.
These fully discrete finite volume methods, described in Section \ref{sec:mol} and \ref{sec:muscl},
are referred to as the base schemes.  Both schemes require linear 
reconstruction on grid cells,
outlined in Section \ref{sec:limit}.  
Finally, both second order schemes evaluate the flux at
the boundary in the examples of Section \ref{sec:compResults} by extrapolating the
pressure to the boundary midpoint. 
We do not discuss domain boundary conditions in this paper,  since 
these procedures are  standard, and do not change due to cut cells.

\subsection{Method of lines} \label{sec:mol}

After generating the cut cell mesh, we approximate the solution to \eqref{eq:conslaw2D} on the base grid using a finite volume method of the form
\begin{equation}\label{eq:fvscheme}
\frac{d}{dt}\mathbf{U}_{i,j} =- \frac{1}{V_{i,j}} \int_{\partial \Omega_{i,j}} \mathbfcal{F} ^* \cdot \mathbf{n} ~dl,
\end{equation}
where $\mathbf{U}_{i,j}$ is a vector of the cell averages on $\Omega_{i,j}$, $\partial \Omega_{i,j}$ is the cell boundary, $\mathbf{n}$ is an outward facing normal, and $\mathbfcal{F}^*$ is a numerical flux function.

Second order accuracy in space is achieved by reconstructing a gradient on each
cell and using it to evaluate the numerical flux at the face midpoints,
illustrated in Figure \ref{fig:2dfig}.  
In our numerical experiments, we use the local Lax-Friedrichs numerical 
flux to solve the Riemann problem at cell interfaces. The integral in \eqref{eq:fvscheme} is  approximated using the midpoint rule.

Second order accuracy in time is obtained by integrating \eqref{eq:fvscheme} using Heun's method.  This is a two-stage Runge Kutta method that can be written
\begin{equation}\label{eq:molscheme}
\begin{aligned}
	\mathbf{U}^{(1)} &= \mathbf{U}^{n} + \Delta t L(\mathbf{U}^n), \\
	\mathbf{U}^{(2)} &= \mathbf{U}^{(1)} + \Delta t L(\mathbf{U}^{(1)}), \\
	\mathbf{U}^{n+1} &= \frac{1}{2}( \mathbf{U}^{n} + \mathbf{U}^{(2)} ) ,	
\end{aligned}
\end{equation}
where $\mathbf{U}^{n}$ is the vector of solution averages on the entire cut cell mesh at time $t^n$,
$\mathbf{U}^{(1)}$,$\mathbf{U}^{(2)}$, are intermediate stages, and $L$ is the operator that results
from discretizing the right-hand-side of \eqref{eq:fvscheme}.
We derive the maximum stable time step on the base grid by using \eqref{eq:molscheme} to solve the linear advection equation on a Cartesian grid with advection velocity $(a,b)$.
Numerical evaluation of the amplification factor that results from a linear stability analysis reveals that a stable time step satisfies
\begin{equation}\label{eq:vn1}
\Delta t   \left( \frac{|a|}{\Delta x} + \frac{|b|}{\Delta y} \right)\leq 1.
\end{equation}
For a discrete maximum principle to be satisfied, a tighter 
time step restriction 
\begin{equation}
\Delta t  \left( \frac{2\Delta x + 2 \Delta y}{\Delta x \Delta y} \right) \sqrt{a^2 + b^2}\leq 1 .
\end{equation}
is required \cite{giuliani2018analysis}, in addition to a limiter to prevent new extrema.
For the Euler equations, the advective speeds $a$ and $b$ are replaced
by the largest magnitude $x$ and $y$ velocities  plus sound speed.

For this method, we apply SRD stabilization to the stage updates
$U^{(1)}$ and $U^{(2)}$ in \eqref{eq:molscheme}, before they are added to 
form $U^{n+1}$.  This is needed since the tiniest cut cells may have
non-physical quantities that need to be stabilized before  the
next stage.  For example, the first step in approximating the Euler
equations is to convert the conserved variables to primitives variables.
This would break down for negative densities. SRD however does not break
down and conservatively and accurately adjusts these values.

\subsection{MUSCL scheme} \label{sec:muscl}
The MUSCL scheme is a one step method that is second order accurate in
space and time. A series of MUSCL schemes  was originated by van Leer 
\cite{vanleer:muscl}. The version we use\footnote{Thanks to Phil 
Colella for the original Cartesian mesh code and Riemann solver,  and for helpful discussions on the shear
layer instability and artificial viscosity.}
is due to Colella \cite{Colella:Unsplit}.
The method is briefly sketched here so that we can describe how it was
adapted for cut cells. 

On a regular cell $(i,j)$ the interface values on the 
faces are computed at the half-step in time $t^{n+1/2}$, and the left and right states
are passed to a Riemann
solver to compute the fluxes.
Using a Taylor series in space and time to second order,
and using the conservation law \eqref{eq:conslaw2D}  to replace the derivative in time, gives the value for the right cell
interface $(i+1/2,j)$ 
\begin{subequations}
\begin{align}
\label{eqn:taylor}
\mathbf{U}_{i+1/2,j}^{n+1/2}    
              &= \mathbf{U}_{i,j}^n +  
\frac{\Delta t}{2} \frac{\partial \mathbf{U}_{i,j}^n}{\partial t} + 
\frac{\Delta x}{2} \frac{\partial \mathbf{U}_{i,j}^n}{\partial x} \\[.08in]
              &= \mathbf{U}_{i,j}^n - \frac{\Delta t}{2} \, 
             \frac{\partial \mathbf{G}_{i,j}^n}{\partial y}  -
            \left( \frac{\Delta t}{2} 
            \frac{\partial \mathbf{F}_{i,j}^n}{\partial \mathbf{U}^n_{i,j}} -
             \frac{\Delta x}{2} \right) \,\frac{\partial \mathbf{U}_{i,j}^n}{\partial x}, 
\end{align}
\end{subequations}
The values at the other interfaces of full cells are similarly defined.

In the original method,
Riemann problems are solved in the transverse direction, e.g. between
centroids $(i,j)$ and $(i,j-1)$, at the edge $j-1/2$, to produce
$\mathbf{G}_{i,j-1/2}$. 
The terms were then differenced to compute 
\begin{equation}
\partial \mathbf{G}_{i,j}/\partial y =  (\mathbf{G}_{i,j+1/2} - \mathbf{G}_{i,j-1/2})/\Delta y .
\label{eqn:transdiff}
\end{equation}
On a uniform mesh 
the first order accurate errors in computing \eqref{eqn:transdiff} 
cancel, and the term itself is 
multiplied by $\Delta t$ in \eqref{eqn:taylor}.
For more details the reader is referred to \cite{Colella:Unsplit}.

At the cut cells the above procedure is no longer accurate,
since  the centroid values are not coordinate aligned near
the cut cells.
We make two modifications to the
computation of $\partial \mathbf{G}/\partial y$, known as the
transverse derivative since it is in the vertical direction
when computing the flux $\mathbf{F}$ in the horizontal direction.
First,
the solution is reconstructed in the transverse direction
to the edge midpoint so it is properly centered in cells that
are adjacent to a cut cell.  This is the situation
in Figure \ref{fig:2dfig} for cell $(i,j+1)$, for example.

Second, many cut cells will not have both edges in the
transverse ($y$) direction. Instead, for all cut cells we
instead compute
$ \partial \mathbf{G}_{i,j}^n / \partial y = ( \partial \mathbf{G}_{i,j}^n / \partial \mathbf{U}_{i,j}^n)( \partial \mathbf{U}_{i,j}^n / \partial y)$,
like the horizontal fluxes in \eqref{eqn:taylor}. This is linearly exact in the cut cells, 
if the gradients themselves are.
We also experimented with dropping this term in the cut cells
altogether. There was no stability problem with this, but it
does introduce an unnecessary difference from
the interior scheme, and would not be linearly exact.

\textit{Note:} it is the transverse derivative term that provides the so-called corner 
coupling, i.e. inclusion
of corner cells in the stencil. This is what  gives  the MUSCL scheme a linear stability limit of
\begin{equation}
\label{eqn:bigcfllimit}
\Delta t \, \max \left (\frac{|a|}{\Delta x} , \frac{|b|}{\Delta y} \right) \leq 1,
\end{equation}
where $(a,b)$ is the advection velocity.  

The trickiest term to adapt to cut cells was an artificial viscosity in the original method that was 
added to each flux, with a
coefficient proportional to the negative divergence of the flow.  The original
code used a large stencil to compute this divergence. We instead use a centered
difference to compute $u_x$ and $v_y$, and where possible, and  
take the max of this over a $3 \times  3$
neighborhood centered around each cell, so that cut cells get this dissipation too. 

\commentout{
As remarked above, it is the transverse derivatives that allows for the larger time
step given in \eqref{eqn:bigcfllimit}.
However, in the neighborhood of a shock the derivatives will be limited, and
possibly set to zero.
If these terms were not used  in the volume mesh, the time step would be reduced to 
\begin{equation}
\Delta t \, \left (\frac{u+c}{\Delta x} + \frac{v+c}{\Delta y} \right) < 1
\end{equation}
which could be as small as half the larger limit in eq. \eqref{eqn:bigcfllimit}.
However, as shown in \cite{mjb:stability2} for one space dimension, 
boundary cells can have
a local {\em cfl} number that is up to twice the stable {\em cfl} of the regular
mesh and the overall scheme remains stable.  We have not found any stability
problems due to limiting of this term.  
}

The multi-dimensional MUSCL scheme due to Colella has several additional
features to robustly handle strong shocks, such as not including terms in
predicting the interface state from characteristics that propagate 
away from the interface. These
steps do not change at the cut cells, so are not discussed here.  

Since MUSCL is a one-step scheme, the SRD stabilization is applied directly
before the final update
$\mathbf{U}^{n+1} = SRD(\mathbf{U}^{n} + \Delta t
L(\mathbf{U}^{n}))$, where $L$ is now the MUSCL operator.

\subsection{Gradient reconstruction and limiting }\label{sec:limit}

The computation of gradients, and for problems with discontinuities limiting
those gradients, arises independently of the finite volume scheme used. 
On the domain interior, when the stencil is regular and does not contain cut cells, standard schemes can be used.
In all our examples we use monotonized central (MC) differencing in both $x$ and $y$ directions.  The MC limited slope in the $x$ direction is
\begin{equation}
\sigma^n_{x,i,j} =  \begin{cases} 
\min \left ( \,  \lvert{ D_c}\rvert,\,
2 \lvert {D_+}\rvert,\,
2 \lvert{D_-}\rvert \,  \right ) \,\times 
\text{ sign } D_c, \quad \text{if} \;\;  D_+ D_- >  0,\\
0 \hspace*{2.8in} \text{otherwise}.
\end{cases}
\end{equation}
Here we use the standard differencing notation
$D_c = (U^n_{i+1,j}-U^n_{i-1,j})/(2 \, \Delta x)$ for the second order accurate central difference and
$D_+ = (U^n_{i+1,j}-U^n_{i,j})/\Delta x$,
$D_- = (U^n_{i,j}-U^n_{i-1,j})/\Delta x$ for the one-sided differences.  The MC limited slope in the $y$ direction is similarly defined.

For cut cells 
we use a least squares gradient reconstruction algorithm, a standard procedure
for unstructured meshes.
A linear reconstruction of the solution on these cells is of the form
\begin{equation}
u^n_{i,j}(x,y) = U_{i,j}^n + \sigma^n_{x,i,j} \,(x-x_{i,j}) +
                     \sigma^n_{y,i,j}\,(y-y_{i,j}),
\label{eqn:lls}
\end{equation}
where $(i,j)$ is the index of either a cut cell or cell with an
irregular stencil, $(\sigma^n_{x},\sigma^n_{y})_{i,j}$ and $(x,y)_{i,j}$
are its gradient and cell centroid, respectively. The least squares
procedure finds the gradient that minimizes the $L_2$ residual when
evaluating  $u^n_{i,j}(x,y)$ at  neighboring cell centroids. 

In this work, we consider both first and second order accurate 
gradients.  The reconstructed first order gradient satisfies in the 
least squares sense
\begin{equation}\label{eqn:linrecon_base}
\sigma^n_{x,i,j}(x_{r,s} - x_{i,j}) +
\sigma^n_{y,i,j}(y_{r,s} - y_{i,j})=
U^n_{r,s} - U^n_{i, j} \quad \forall (r,s) \in R_{i,j},
\end{equation}
where $R_{i,j}$ is the set of cell indices used for slope reconstruction on cell $(i,j)$ in the 
base scheme.  Here, $R_{i,j}$ is the $3\times 3$ neighborhood  centered on $(i,j)$.
The reconstructed second order gradient satisfies in the least squares sense
\begin{equation}
\begin{aligned}\label{eqn:linrecon_base2}
&\sigma^n_{x,i,j}(x_{r,s} - x_{i,j}) +
\sigma^n_{y,i,j}(y_{r,s} - y_{i,j})  + \\
&\quad \frac{1}{2}\sigma^n_{xx,i,j}[(x_{r,s} - x_{i,j})^2 - S_{xx,i,j}]  + \\
& \quad \; \sigma^n_{xy,i,j}[(x_{r,s} - x_{i,j})(y_{r,s} - y_{i,j})-S_{xy,i,j}] +\\
&\quad  \frac{1}{2}\sigma^n_{yy,i,j}[(y_{r,s} - y_{i,j})^2 - S_{yy,i,j}]
 \; = \;  U^n_{r,s} - U^n_{i, j} \quad \forall (r,s) \in R_{i,j},
\end{aligned}
\end{equation}
where $\sigma^n_{xx,i,j}$, $\sigma^n_{xy,i,j}$, $\sigma^n_{yy,i,j}$ 
are quadratic degrees of freedom and are discarded.  In this case, $R_{i,j}$ is either the $3\times 3$ tile centered on $(i,j)$ when $(i,j)$ is a whole cell neighboring a cut cell, or the $5\times 5$ tile centered on $(i,j)$ when $(i,j)$ is a cut cell.
Note that a cut cell needs a larger neighborhood because 
approximately half of its cells are not in the flow domain.  

Regular cells that are adjacent to a cut cell will also need special treatment to
compute a second order accurate gradient. We have experimented with three  approaches
and found almost indistinguishable results. The simplest is to use the 
procedure mentioned above for cut cells
- fit a least squares polynomial in the $3 \times 3$ neighborhood centered on that cell.
We also tried $5 \times 5$ neighborhoods, in hope of smoother transitions between cut cells
and the interior cells.  Finally, we tried using a centered
gradient in only one dimension, if there was one, and using a recentering approach to 
compute the second order accurate difference in the other direction.  This has an overall
smaller stencil, but still uses the $3 \times 3$ neighborhood
if there is no regular direction,
and involves more testing.  

For problems with discontinuities, the gradient will need to be limited
to prevent overshoots and retain positivity for quantities like density and
pressure.
We use the Barth Jespersen (BJ)  limiter \cite{barth-jespersen} to limit on 
cut cell grids. 
This is a scalar limiter, where both $\sigma^n_{x,i,j}$ and $\sigma^n_{y,i,j}$ 
are reduced by the same scalar to prevent new extrema.  
We compute the minimum and maximum values over the reconstruction 
stencil $R_{i,j}$, 
\begin{equation} 
m_{i,j} = \max_{(r,s) \in R_{i,j}} U^n_{r,s} \text{ and } 
M_{i,j} = \max_{(r,s) \in R_{i,j}} U^n_{r,s}.
\label{eqn:bj1}
\end{equation}
The reconstructed gradient on cell $(i,j)$ is limited by a non-negative 
scalar $\alpha \in [0,1]$, so that when ${u}_{i,j}(x,y)$ 
is evaluated at the centroids of the neighborhoods in $R_{i,j}$ it
lies between $m_{i,j}$ and $M_{i,j}$.
(Apologies for reusing the symbol $\alpha$, since it is commonly 
used to describe Barth-Jespersen-type limiters, as well as the 
mesh width $\alpha h$  of small cells).

The limited numerical solution is
\begin{equation}
     \tilde{u}^n_{i,j}(x,y) = U_{i,j}^n + \alpha \, [{\sigma}^n_{x,i,j} ( x -  x_{i,j}) \, 
   + {\sigma}^n_{y,i,j}( y -  y_{i,j})].
\end{equation}
Define
\begin{equation}\label{eq:bj_alpha}
    \alpha_{r,s} = \begin{cases}
           \min \left(1,\frac{M_{i,j}-U_{i,j}^n}{U^n_{r,s} - U_{i,j}^n} \right)
    \quad  \text{ if } \,   U_{r,s}^n - U_{i,j}^n >  0,\\[.08in]
            \min \left(1, \frac{m_{i,j}-U_{i,j}^n}{U^n_{r,s} - U_{i,j}^n} \right)  
    \quad  \text{ if }  \, U^n_{r,s} - U_{i,j}^n < 0.\\[.08in]
             1    \hspace*{1.45in}  \text{if} \; \, U^n_{r,s} - U_{i,j}^n = 0.
    \end{cases}
\end{equation}
Then choose
\begin{equation}\label{eqn:alpha}
\alpha = \min_{(r,s) \in R_{i,j}} \alpha_{r,s} .
\end{equation}

By reconstructing to the neighboring cell centroids, this procedure is linearity preserving.  

The approach described above differs slightly from the original Barth-Jespersen limiter in
\cite{barth-jespersen}, where the numerical solution is reconstructed to 
points on cell interfaces.
This is often the procedure used on unstructured meshes. However cut cell meshes are much more
irregular and without this fix BJ can lead to much less accurate solutions.

\section{State redistribution in two dimensions}\label{sec:srdAlg}

This section describes the second order accurate
state redistribution algorithm. 
We will show that SRD preserves linear functions and is conservative, since
it is not obvious that the unusual weightings in our algorithm
preserve these important properties.

\subsection{State redistribution preprocessing}\label{sec:preprocessing}

In this section, we describe mesh dependent quantities that will be used when 
applying SRD. 
Each cut cell needs two pieces of information: the 
cells that belong to its own merging neighborhood,
and how many neighborhoods it belongs to.
This is the two-dimensional analogue of the one-dimensional 
nonuniform grid preprocessing presented in Section \ref{sec:srd1d}, now done on cut cell grids.
The preprocessing determines these quantities: merging neighborhoods and overlap counts, 
weighted volumes and centroids.
For moving geometries, they would need to be be recomputed
when the geometry is modified.

\subsubsection*{Merging neighborhoods}

In one dimension, neighborhoods can be defined by merging a small cell with neighbors on its left, on its right, or both.  Of these three approaches, the last was used on the nonuniform grid presented in Section \ref{sec:srd1d}.
In two dimensions, there is much more freedom in defining a merging neighborhood.  
We investigated two different ways: (1) normal merging and (2) centered merging.

\begin{figure}[h]
    \centering
    \includegraphics[width=0.5\linewidth]{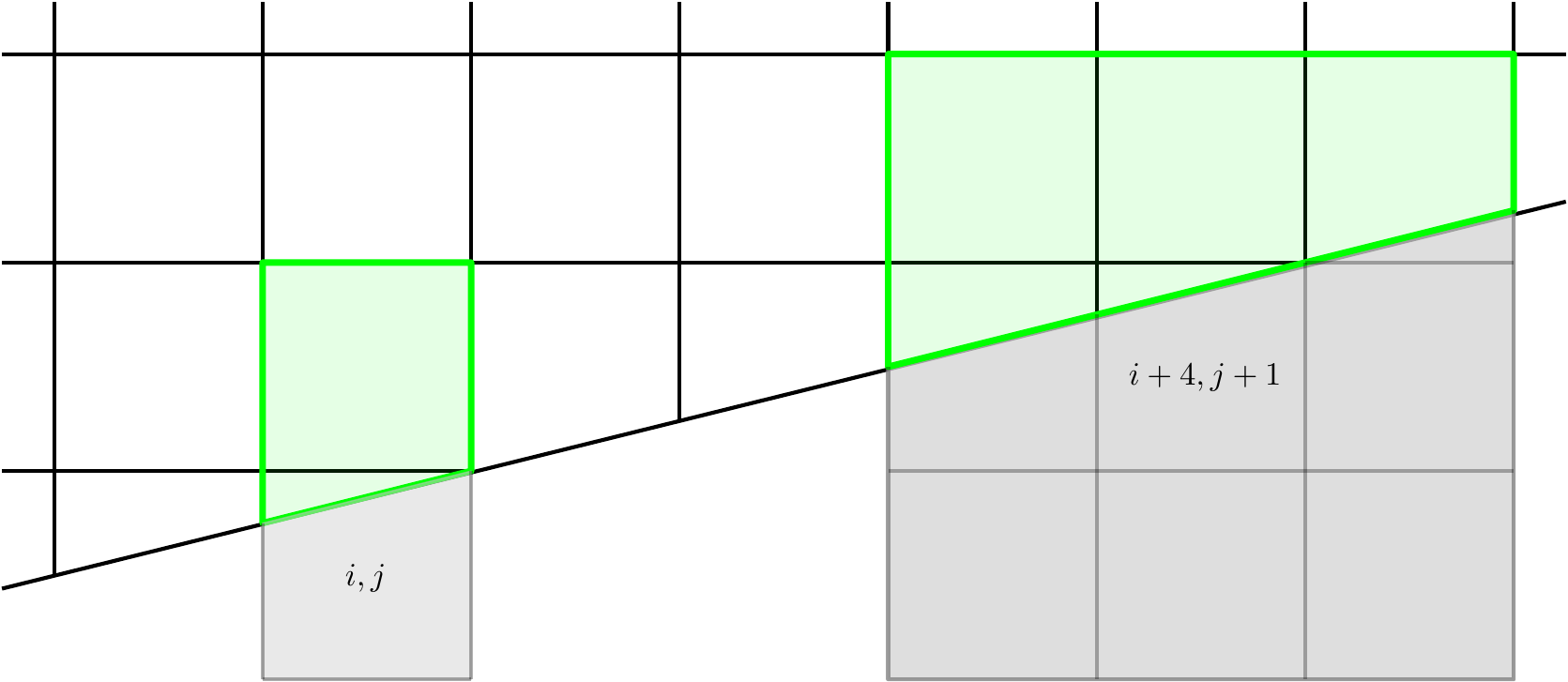}
    \caption{\sf On the left, a small cell is merged with a cell in the direction closest to the 
    	boundary normal.  On the right, a small cell is merged with neighbors that are at most one cell away, i.e., cells located in the $3\times3$ tile centered on $(i+4,j+1)$.}
    \label{fig:neighborhoods}
\end{figure}

Normal merging associates cut cells with neighbors in the direction closest to the 
boundary normal.  This is illustrated in Figure \ref{fig:neighborhoods},
where the normal merging neighborhood associated with cut cell $(i,j)$ is highlighted in green.

Centered merging associates cut cells with neighbors that are 
symmetrically located in each direction around the cell.
These cells are located in the $3 \times 3$ tile centered on the small cell.  Figure \ref{fig:neighborhoods}
illustrates this for cut cell $(i+4,j+1)$, where the centered neighborhood is highlighted in green.

The larger the neighborhood the more diffusive the results. Thus we do not want neighborhoods that 
are too large.  However, they must satisfy a size constraint to dampen unstable growth in the numerical 
solution. To this end, we determine neighborhoods, either normal or centered, such that the
volume of the neighborhood is at least half the volume  of an uncut cell, i.e., 
\begin{equation} \label{eqn:vmerge}
\sum_{(r,s) \in M_{i,j}} V_{r,s} \geq \frac{1}{2}\Delta x \Delta y,
\end{equation}
where $M_{i,j}$ denotes the set of cell indices that belong to 
merging neighborhood $(i,j)$.    
In one dimension it has been shown \cite{mjb:stability2} that a cut cell at the boundary 
that is at least half the regular cell size  is stable using a full time step $\Delta t$.  
We have not encountered any issues with this choice in two dimensions.

\begin{figure}[h]
\hspace*{.5in}
	\subfloat[\sf Normal neighborhood (in red) for the cut cell $(i,j)$.]{\includegraphics[width=.30\textwidth]{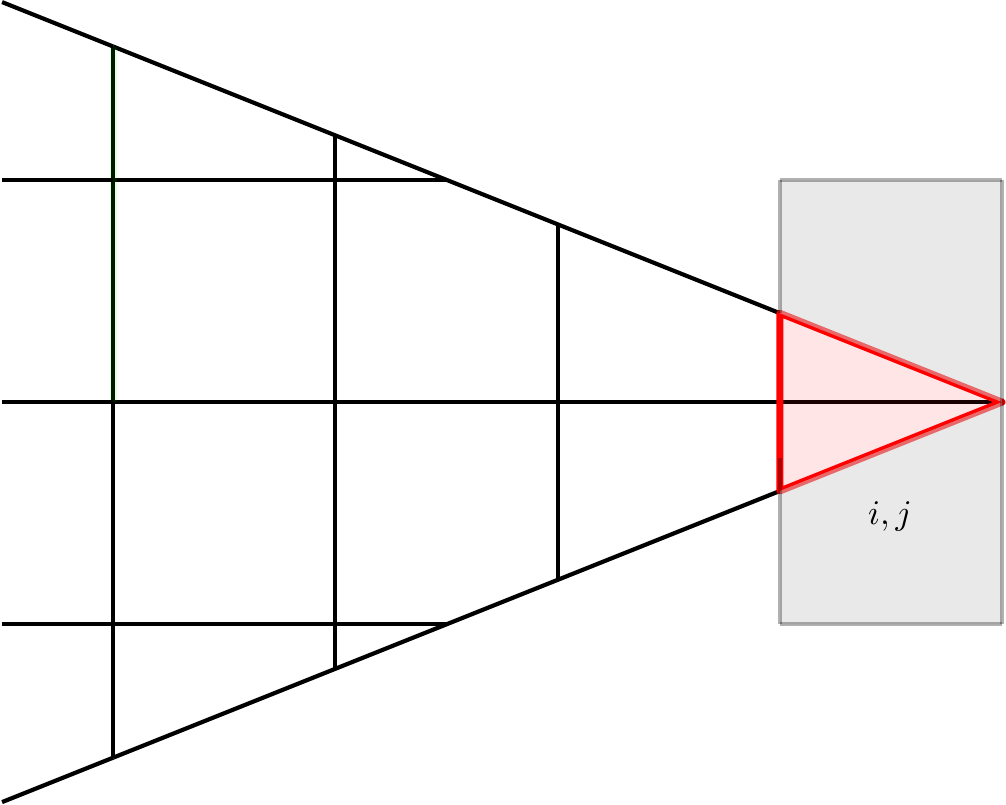} \label{fig:normalneighborhood}}
	\hfill
	\subfloat[\sf $3\times 3$ merging neighborhood (in green) for centered on cut cell $(i,j)$.]{\includegraphics[width=.35\textwidth]{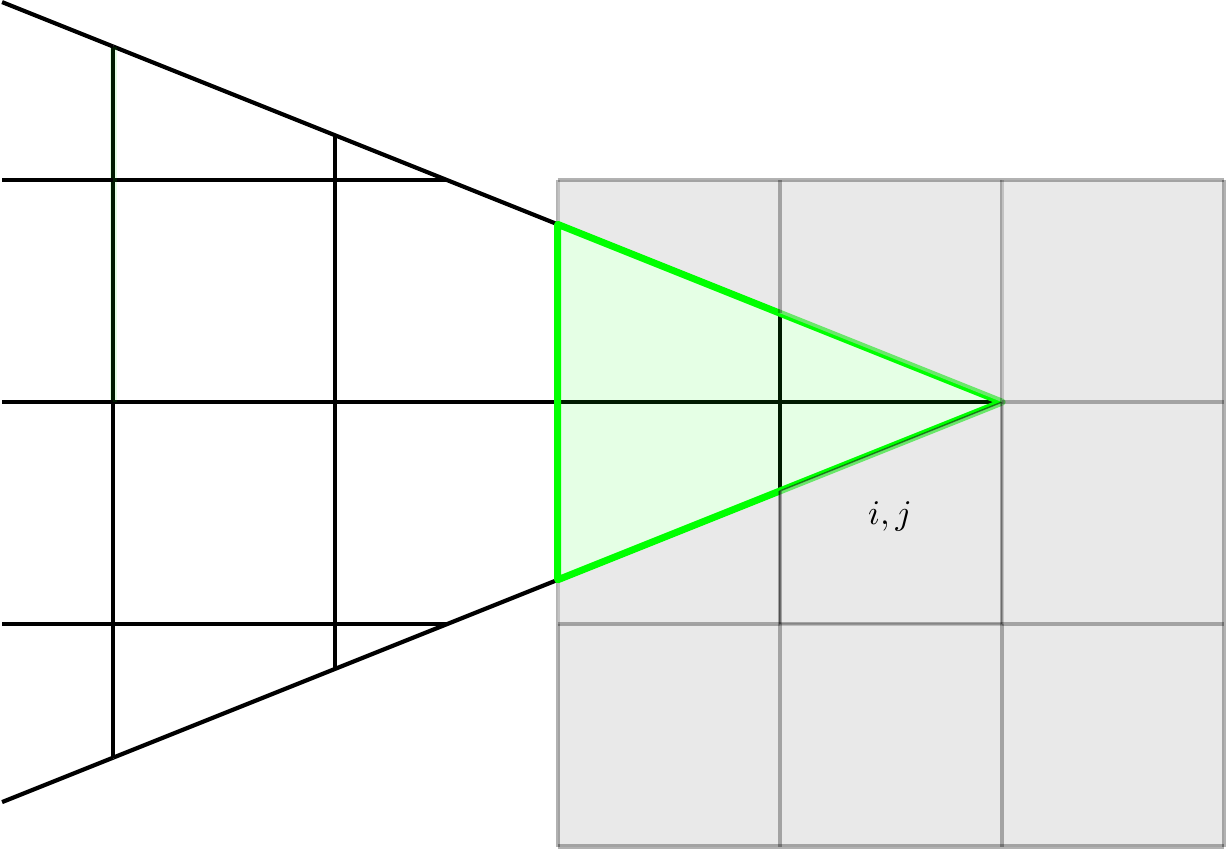} \label{fig:3x3neighborhood}}
	\caption{\sf It can happen that merging only in the normal neighborhood is not 
        large enough, e.g. if the small cell merges with another small cell, and \eqref{eqn:vmerge} is not satisfied.  In this case, we use the $3\times 3$ neighborhood. }
\end{figure}

For smooth solutions a centered neighborhood is sufficient, but for shocks this yields unsatisfactory results. Therefore, we use the normal neighborhood everywhere possible.
The normal neighborhood cannot be used when, e.g., a neighboring cell is also cut and the
merging neighborhood is not sufficiently large (Figure \ref{fig:normalneighborhood}).  In this case, we must use central merging with cells on a $3\times3$ tile (Figure \ref{fig:3x3neighborhood}), or, if that merging neighborhood is not large enough, with cells on the $5 \times 5$ tile.
In this Figure, the gray cells are exterior to the fluid domain, and
are drawn for context.

After forming the merging neighborhoods, each cell counts the number of neighborhoods it belongs to.
For cell $(i,j)$, this neighborhood count is called $N_{i,j}$.

\begin{figure}[h]
	\centering
	\subfloat[Merging neighborhood of cell $(i,j)$.]{\includegraphics[width=0.30\textwidth]{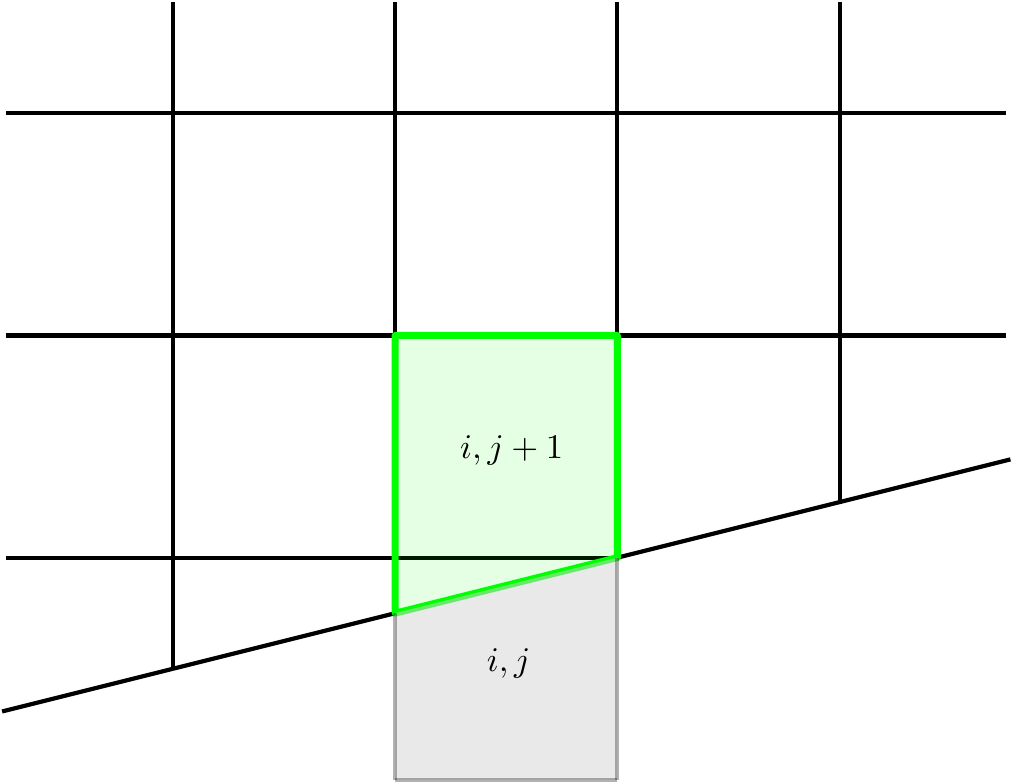} \label{fig:2nborTile1}}
	\qquad
	\subfloat[Merging neighborhood of cell $(i,j+1)$.]{\includegraphics[width=.30\textwidth]{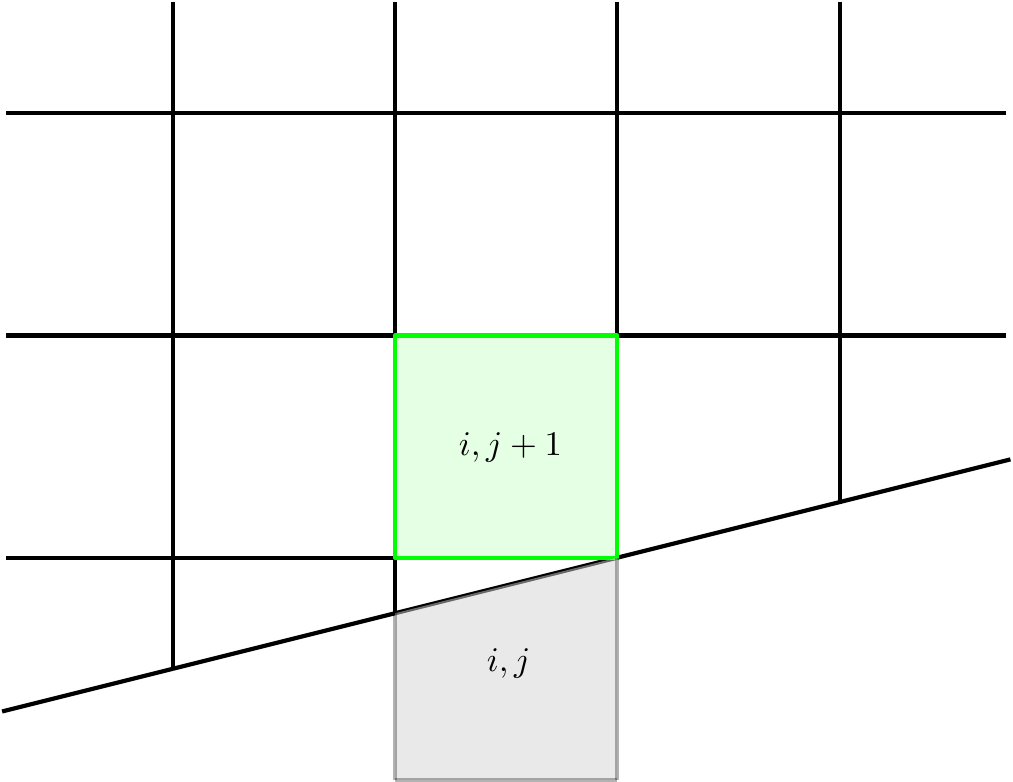} \label{fig:2nborTile2}}
	\caption{\sf Example of two neighborhoods, highlighted in green, that overlap on cell $(i,j+1)$.} \label{fig:2nborTile}
\end{figure}

To more clearly illustrate overlapping merging neighborhoods in two dimensions, consider the cut cell 
mesh in Figure \ref{fig:2nborTile}.   The neighborhoods $(i,j)$ and $(i,j+1)$ overlap on cell $(i,j+1)$ in 
the base grid.  The normal merging neighborhood of cut cell $(i,j)$ is highlighted in green in Figure \ref{fig:2nborTile1}.  
Since cell $(i,j)$ does not satisfy the volume constraint in \eqref{eqn:vmerge},
the neighborhood of $(i,j)$, $M_{i,j}$, must include both $(i,j)$ and $(i,j+1)$.  
The neighborhood associated with $(i,j+1)$ in highlighted in green in Figure \ref{fig:2nborTile2}.  
Since the volume of cell $(i,j+1)$ satisfies the volume constraint in \eqref{eqn:vmerge}, $M_{i,j+1}$, only contains cell $(i,j+1)$.
  
\subsubsection*{Weighted centroids and volumes}

The neighborhood's weighted volume is defined as
\begin{equation}
\label{eqn:voldef}
{\widehat V}_{i,j} =  \sum_{(r,s) \in M_{i,j} } \,  \frac{V_{r,s}}{N_{r,s}}.
\end{equation}
where $M_{i,j}$ is the set of cells in the neighborhoods associated with cell $(i,j)$. 
The weighted centroid of the merging neighborhood is defined as
\begin{equation}
\label{eqn:centroiddef}
({\widehat x}_{i,j},{\widehat y}_{i,j}) = \frac{1}{\widehat V_{i,j}} \sum_{(r,s) \in M_{i,j} } \,  \frac{V_{r,s}}{N_{r,s}}(x_{r,s},y_{r,s}),
\end{equation}

In general, the weighted volume is not the physical volume of 
the merging neighborhood 
unless the neighborhood  is not overlapped by any other neighborhoods 
(for example, cells $-3$ and $3$ in Figure 2).  Similarly, the weighted and physical 
centroids are different unless the overlap count on each cell in the neighborhood is the 
same, $N_{r,s} = N_{p,q} ~ \forall (r,s), (p,q) \in M_{i,j} $.  
This is because the overlap counts in \eqref{eqn:centroiddef} would cancel, leading to identical  weighted and physical 
centroids of the neighborhood.




\subsection{State redistribution postprocessing} \label{sec:srd_postprocessing}

In this section, we describe postprocessing on two-dimensional meshes with overlapping neighborhoods.
This is the analogue of postprocessing on one dimensional nonuniform grids in Section \ref{sec:srd1d}.
The SRD stabilization is applied after each stage for the method of lines (Section \ref{sec:mol}) 
or time step for the MUSCL scheme (Section \ref{sec:muscl}), denoted generically as
\begin{equation} \label{eq:stage_step}
\widehat{U} = U^n + \Delta t  L(U^n).
\end{equation}
We refer to
$\widehat{U}$ as the provisionally updated solution.

\subsubsection*{Step 1. Compute weighted solution averages on each neighborhood}   

The solution average on each neighborhood is given by
\begin{equation}
\label{tiledef}
\widehat{Q}_{i,j} =  \frac{1}{{\widehat V}_{i,j}} \, \sum_{(r,s) \in M_{i,j}} \,  
\frac{V_{r,s}}{N_{r,s}}  \widehat{U}_{r,s}
\end{equation}
where the volume ${\widehat V}_{i,j}$ is the weighted neighborhood volume defined in \eqref{eqn:voldef}
and $N_{r,s}$ is the number of neighborhoods  associated with cell $i,j$ defined in 
Section \ref{sec:preprocessing}.  
Analogous to the one dimensional case in Section \ref{sec:srd1d}, the weighted solution 
averages are a convex combination of provisional solution averages on the
neighborhood associated with cell $(i,j)$.  

\subsubsection*{Step 2. Reconstruct and limit a gradient on each neighborhood}

For second order accuracy in space we need to compute a gradient on each 
neighborhood, again  using a least squares procedure.
The set of merging neighborhood indices for reconstruction on the
neighborhood associated with cell $(i,j)$ is called $\widehat R_{i,j}$.
We note that this set need not be the same as $R_{i,j}$ used in the cut cell gradient reconstruction on the base grid (Section \ref{sec:limit}). 
Similar to the base grid reconstruction in \eqref{eqn:lls}, the reconstruction on neighborhood $(i,j)$ is of the form
\begin{equation}\label{eq:qrecon}
\widehat{q}_{i,j}(x,y) = \widehat{Q}_{i, j} + \widehat{\sigma}_{x,i,j}(x - \widehat{x}_{i,j}) + \widehat{\sigma}_{y,i,j}(y - \widehat{y}_{i,j}),
\end{equation}
where $\widehat{Q}_{i, j}$ is the weighted neighborhood average defined in \eqref{tiledef}, $(\widehat{x}_{i,j},\widehat{y}_{i,j})$ is the weighted centroid of neighborhood $(i,j)$, and $(\widehat{\sigma}_{x,i,j},\widehat{\sigma}_{y,i,j})$ is 
the gradient on the merging neighborhood.
The neighborhood gradient $(\widehat{\sigma}_{x,i,j},\widehat{\sigma}_{y,i,j})$ satisfies in the least squares sense
\begin{equation}\label{eqn:linrecon}
\widehat{\sigma}_{x,i,j}(\widehat{x}_{r,s} - \widehat{x}_{i,j}) +
\widehat{\sigma}_{y,i,j}(\widehat{y}_{r,s} - \widehat{y}_{i,j})=
\widehat{Q}_{r,s} - \widehat{Q}_{i, j} \quad \forall (r,s) \in \widehat{R}_{i,j}.
\end{equation}
Note that \eqref{eqn:linrecon} uses the neighborhood's weighted centroids $(\widehat{x}_{r,s},\widehat{y}_{r,s})$ instead of the cell centroids.  In Section \ref{sec:linex}, we prove that this procedure is linearity preserving.
Alternatively, second order gradients on the neighborhoods can be obtained by fitting a quadratic
as in Section \ref{sec:limit}, and discarding the second derivative terms.
We will compare these procedures in the computational results in Section \ref{sec:ssv}.

\begin{figure}
    \centering
    \subfloat[The blue merging neighborhood is associated with cut cell
    $(i,j)$.  The green merging neighborhoods are associated with cut cells $(i-1,j)$ and $(i+1,j)$.]{\includegraphics[width=0.30\linewidth]{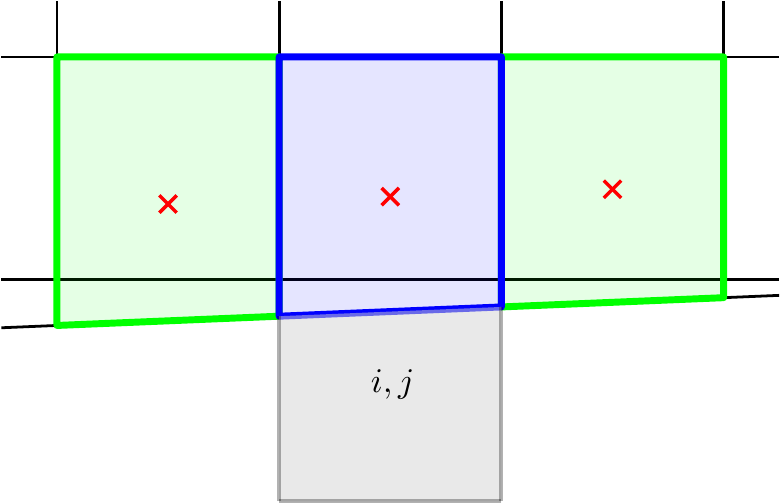} \label{fig:a}} \hfill
    \subfloat[The green merging neighborhoods are associated with the whole cells $(i-1,j+1)$, $(i,j+1)$, and $(i+1,j+1)$.]{\includegraphics[width=0.30\linewidth]{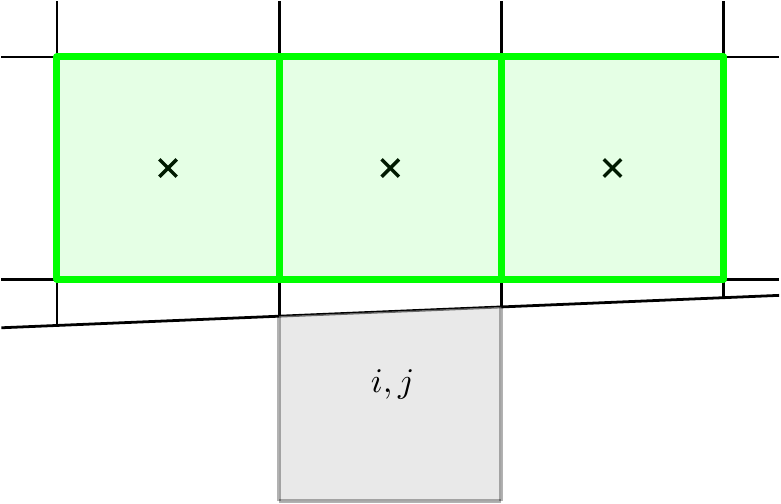} \label{fig:b}} \hfill
    \subfloat[The weighted centroids of the reconstruction neighborhoods in Figures \ref{fig:a}, \ref{fig:b}.]{\includegraphics[width=0.30\linewidth]{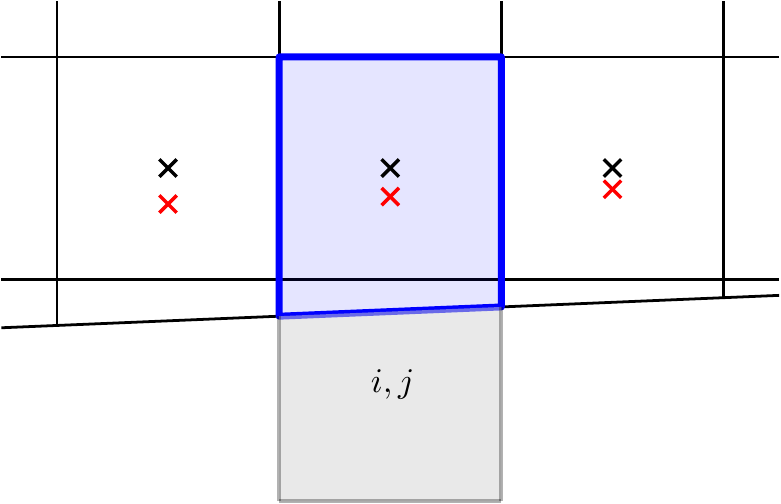} \label{fig:c}} 
    \caption{\sf 
    The weighted centroids of the neighborhoods in $\widehat{R}_{i,j}$ are indicated with a cross ($\times$).  
    Here $\widehat{R}_{i,j}$ is the set of merging neighborhoods associated
    with cells on the $3\times3$ tile centered on cell $(i,j)$.
    Clearly, the weighted centroids are too close to one another in the $y$ direction (Figure \ref{fig:c}).  It follows that least squares system to reconstruct the $x$ and $y$ derivatives of the numerical solution on the blue merging neighborhood is very ill-conditioned.
}
    \label{fig:tooclose}
\end{figure}

The set of neighborhoods used for gradient reconstruction on neighborhood $(i,j)$, $\widehat R_{i,j}$, is the $3 \times 3$ tile.
It could happen that this set does not contain enough neighborhoods, or that the weighted centroids of these neighborhoods are too close to compute a well-conditioned gradient.
This is the case in Figure \ref{fig:tooclose}, where the weighted centroids are too close in the $y$ direction.

We remedy this by increasing the stencil size for the gradient computation 
if the neighborhood does not contain another weighted centroid at least  
$\frac{1}{2}\Delta x$  and $\frac{1}{2}\Delta y$ away  in the $x$ or $y$ 
direction respectively.
For example, if the weighted centroids are too close in the $x$ 
direction, but not the $y$ direction, then the $5\times 3$ tile is used as the 
reconstruction neighborhood.  Similarly, if the weighted centroids are too 
close in the $y$ direction, but not the $x$ direction, then the $3\times 5$ tile is 
used as the reconstruction neighborhood.  The neighborhood size is increased until this distance requirement is satisfied in both $x$ and $y$ directions.  In Figure \ref{fig:tooclose},
the appropriate reconstruction neighborhood is the $3\times 5$ reconstruction tile.

To limit the reconstructed gradient we apply the BJ limiter, this time over $\widehat{R}_{i,j}$
rather than $R_{i,j}$ in \eqref{eqn:bj1} and \eqref{eqn:alpha}.  This procedure is also 
linearity preserving. 

\subsubsection*{Step 3. Final solution update} 

The final update on cell $(i,j)$ is the average of all its neighborhood reconstructions 
evaluated at $(i,j)$'s physical centroid $(x_{i,j},y_{i,j})$. 
This is given by 
\begin{equation} \label{eqn:final_update_linear}
U^{n+1}_{i,j} =   \frac{1}{N_{i,j}}\sum_{(r,s)  \in W_{i,j}}\hat{q}_{r,s}(x_{i,j},y_{i,j}).
\end{equation}
where $W_{i,j}$ is the list of indices of the neighborhoods that overlap cell $(i,j)$.

\vspace*{.5in}


\noindent\textit{Note }: The final update formula \eqref{eqn:final_update_linear} can easily be implemented with a 
nested for loop as in Algorithm \ref{alg:finalupdate},  
instead of computing the set $W_{i,j}$ in \eqref{eqn:final_update_linear}.  
The outer loop iterates over the merging neighborhoods $(i,j)$ and the inner loop iterates over each cell $(r,s)$ in neighborhood $(i,j)$.  Each merging neighborhood $(i,j)$ gives a contribution $ \hat{q}_{i,j}(x_{r,s}, y_{r,s})/N_{r,s} $ to the cells $(r,s)$ that belong to it. 
\begin{algorithm}[H]
		\caption{\sf Final solution update} \label{alg:finalupdate}
	\begin{algorithmic}
	\For{$i,j$}
	\State $U^{n+1}_{i,j} \leftarrow 0$
	\EndFor
	\For{$i,j$}
		\For{$(r,s) \in M_{i,j}$}
			\State $U^{n+1}_{r,s} \leftarrow U^{n+1}_{r,s} + \hat{q}_{i,j}(x_{r,s}, y_{r,s})/N_{r,s} $
		\EndFor
	\EndFor
	\end{algorithmic}
\end{algorithm}

\subsection{Linear exactness} \label{sec:linex}
In this section, we show that second order accurate state redistribution preserves linear functions.
One might wonder about this since the centroids and solution averages are weighted in this unusual way.
In addition, the local truncation error does not imply the order of
accuracy of the scheme  as it does on regular meshes \cite{Kreiss:white2}. 
Here we simply show that a
linear function remains exact after SRD, if the base scheme is linearly
exact.

Consider the grid function of the numerical solution after one time step
or stage, $\widehat{U}$.  Assume that it can be written in terms of a linear function $f(x,y)$ of the $x$ and $y$ coordinates, i.e.,
\begin{equation}
    \label{eqn:uhatlinear}
\widehat{U}_{i,j} = f(x_{i,j},y_{i,j}) = ax_{i,j} + by_{i,j} + c.
\end{equation}
This assumption is valid since both the method of lines (Section \ref{sec:mol}) and the MUSCL scheme (Section \ref{sec:muscl}) are linearity preserving.  
From \eqref{eqn:uhatlinear} and the expression for the average on the merging neighborhood $(i,j)$ in \eqref{tiledef}, we have
\begin{equation}
    \label{eqn:linear1}
\widehat{Q}_{i,j} = \frac{1}{{\widehat V}_{i,j}} \, \sum_{(r,s) \in M_i} \,  
\frac{V_{r,s}}{N_{r,s}}  \,\, (ax_{r,s} + by_{r,s} + c).
\end{equation}
Distributing the summation in \eqref{eqn:linear1}, we have
\begin{equation}\label{eqn:linear2}
\widehat{Q}_{i,j} =  a \left(\frac{1}{{\widehat V}_{i,j}} \, \sum_{(r,s) \in M_{i,j}} \,  
\frac{V_{r,s}}{N_{r,s}} x_{r,s} \right) + b\left(\frac{1}{{\widehat V}_{i,j}} \, \sum_{(r,s) \in M_{i,j}} \,  
\frac{V_{r,s}}{N_{r,s}} y_{r,s} \right) + c\left(\frac{1}{{\widehat V}_{i,j}} \, \sum_{(r,s) \in M_{i,j}} \,
\frac{V_{r,s}}{N_{r,s}}\right) .
\end{equation}
From the definition of the weighted centroid and volume of the merging neighborhood
$(i,j)$ in \eqref{eqn:voldef} and \eqref{eqn:centroiddef}, respectively, \eqref{eqn:linear2} becomes
\begin{equation}\label{eqn:linear3}
\widehat{Q}_{i,j} =  a \widehat{x}_{i,j} + b\widehat{y}_{i,j} + c = f(\widehat{x}_{i,j},\widehat{y}_{i,j}).
\end{equation}
Now, on neighborhood $(i,j)$, we solve the least squares system 
\eqref{eqn:linrecon} to find $\widehat{\sigma}_{x,i,j}$ and
$\widehat{\sigma}_{y,i,j}$, the gradient on the merging neighborhood.  Using
\eqref{eqn:linear3} in \eqref{eqn:linrecon}, and due to the linearity of $f$, the following system
\begin{equation}
\widehat{\sigma}_{x,i,j}(\widehat{x}_{r,s} - \widehat{x}_{i,j}) + \widehat{\sigma}_{y,i,j}(\widehat{y}_{r,s} - \widehat{y}_{i,j})= f(\widehat{x}_{r,s}, \widehat{y}_{r,s}) - f(\widehat{x}_{i,j}, \widehat{y}_{i,j}) \quad \forall (r,s) \in \widehat{R}_{i,j},
\end{equation}
is solved exactly by $\widehat{\sigma}_{x,i,j}=a$ and
$\widehat{\sigma}_{y,i,j}=b$.  In other words, the exact gradient 
of $f$ is reconstructed on merging neighborhood $(i,j)$.  
The reconstructed solution is then
\begin{equation}
    \label{eqn:qrecon1}
    \hat{q}_{i,j}(x,y) = f(\widehat{x}_{i,j},\widehat{y}_{i,j}) + a(x-\widehat{x}_{i,j})+b(y-\widehat{y}_{i,j}) .
\end{equation}
Using \eqref{eqn:linear3}, \eqref{eqn:qrecon1} becomes
\begin{equation}
    \label{eqn:qrecon2}
    \hat{q}_{i,j}(x,y) = f(x,y).
\end{equation}
By \eqref{eqn:final_update_linear}, the final solution update is
\begin{equation} 
U^{n+1}_{i,j} = \frac{1}{N_{i,j}}\sum_{(r,s) \in W_{i,j}}f(x_{i,j},y_{i,j}).
\end{equation}
The function values are exact, and there are $N_{i,j}$ of them, so
after dividing by $N_{i,j}$ we get
\begin{equation} 
U^{n+1}_{i,j} = f(x_{i,j},y_{i,j}),
\end{equation}
which shows that second order accurate state redistribution preserves linear functions.

\subsection{Conservation}\label{sec:cons}
In this section, we show that the total mass of the numerical solution before and after state redistribution does not change.  It follows from the final update in \eqref{eqn:final_update_linear} that the total mass after state redistribution is
\begin{equation}\label{eq:mi}
\sum_{i,j} V_{i,j} U^{n+1}_{i,j} = \sum_{i,j} \sum_{(r,s) \in
M_{i,j}}\frac{1}{N_{r,s}} V_{r,s} \widehat q_{i,j}(x_{r,s},y_{r,s}) ,
\end{equation}
where $\widehat q_{i,j}(x)$ is that neighborhood's polynomial reconstruction defined in \eqref{eq:qrecon}.
Using the definition
of $\widehat{q}_{i,j}(x,y)$, we have
\begin{equation}\label{eq:mi1}
\sum_{i,j} V_{i,j} U^{n+1}_{i,j} = \sum_{i,j} \widehat {Q}_{i,j} 
\widehat {V}_{i,j}.
\end{equation}
Using the definition of the weighted average \eqref{tiledef}, we have
\begin{equation}\label{eq:mi2}
\sum_{i,j} V_{i,j} U^{n+1}_{i,j} = \sum_{i,j} \sum_{(r,s) \in M_{i,j} }\frac{V_{r,s}}{N_{r,s}} \widehat U_{r,s}.
\end{equation}
Since $N_{r,s}$ indicates the number of times cell $(r,s)$ is overlapped by merging neighborhoods, it follows that the $\frac{V_{r,s}}{N_{r,s}} \widehat U_{r,s}$ term is repeated $N_{r,s}$ times in the sum of \eqref{eq:mi2}.  Thus, simplifying \eqref{eq:mi2}, it follows that
\begin{equation} \label{eq:final}
\sum_{i,j} V_{i,j} U^{n+1}_{i,j} = \sum_{i,j} V_{i,j} \widehat U_{i,j}
\end{equation}
and the total mass before and after state redistribution is the same.

\section{Computational Results}\label{sec:compResults}

In this section we show several computational experiments using state
redistribution to solve the Euler
equations. We will also use these
examples to examine 
properties of different gradient choices and base schemes. 

\subsection{Supersonic vortex study}\label{sec:ssv}

\begin{wrapfigure}{R}{0.45\textwidth}
\centering
\vspace*{-.2in}
\includegraphics[height=0.45\textwidth]{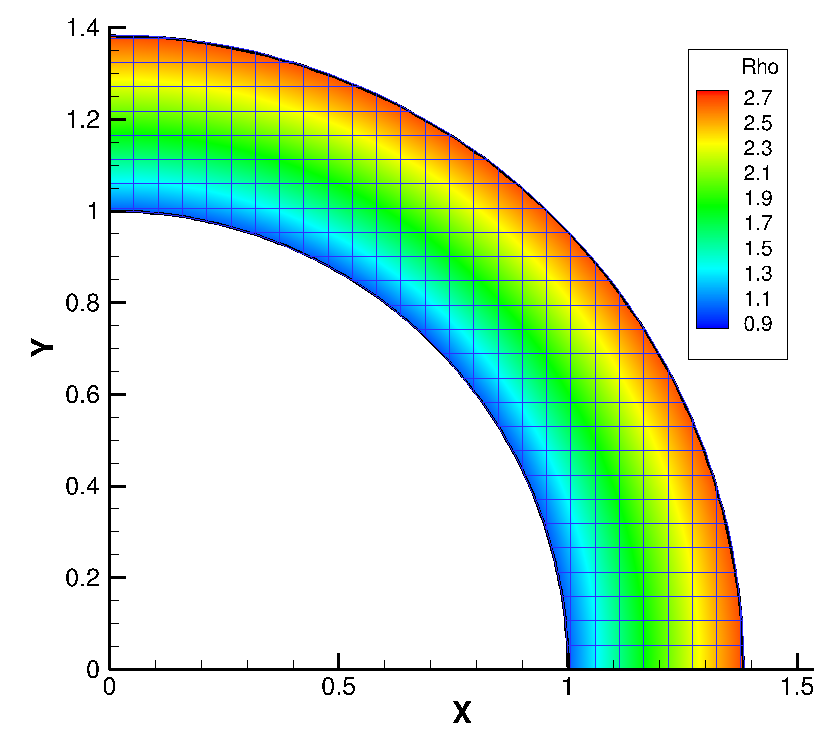}
\caption{\sf The domain is discretized by 27 cells in each direction, and
the exact density profile is shown. Ghost cells and solid cells are not shown.
\label{fig:ssvFig}}
\end{wrapfigure}

We compute the solution to a supersonic flow 
around a quarter circle.  This problem has often been used in accuracy 
studies \cite{aftosmis:acc} since it has an exact solution to the Euler
equations that is smooth, given by 
\begin{equation}
\rho = \rho_i \left \{ 1 + \frac{\gamma-1}{2} \, M_i^2 \left [ 1 - (\frac{r_i}{r})^2
\right ] \right \} ^
{\frac{1}{\gamma-1}}
\end{equation}
and $ u = a_i \, M_i \, (\frac{r_i}{r})\,  \sin (\theta)$, 
$ v = -a_i\,  M_i\,  (\frac{r_i}{r})\,  \cos(\theta)$, and
$ p = \rho^\gamma / \gamma$.
Here, the inner radius is $r_i = 1.0$,  the outer radius
is $r = 1.384$, $\rho_i=1$, and the Mach number on the inner circle
$M_i = 2.25$ in our experiments. 
In this
normalization we use $a_i = 1$, $p_i = 1/\gamma$. 
The second order MOL scheme is used, 
and the time step is chosen using a CFL of 1.0 in \eqref{eq:vn1}
based on the regular cell volume.

We march to steady state, until the maximum density update is 
below $10^{-10}$.  
The solution is smooth, so no limiters are needed.
The exact solution is used to set the ghost cells at the inflow and
outflow boundaries. The domain size is $[0,1.43]$ by $[0.0,1.4301]$ 
(slightly different to prevent mesh  degeneracies).  
The density profile and computational domain is illustrated in Figure
\ref{fig:ssvFig}.

This example will demonstrate the accuracy of SRD.
We also use this example to compare the accuracy of different
formulations for the cut cell and SRD neighborhood  gradients, since
as noted previously, gradients are an important part of this algorithm.
Table \ref{tab:ssv} we compares the accuracy of three such formulations.
For first order accurate gradients, we use a linear least squares 
reconstruction for both the irregular cell gradient (cut cells and their one-away neighbors) 
described in Section \ref{sec:limit}, and the SRD gradients, which 
update the cut cell solution after stabilization. Second order accurate 
gradients fit a quadratic for both
the cut cells and tiles as described in Section \ref{sec:limit}, 
but only the first derivative terms are used. 
As an intermediate experiment, we  fit a quadratic using least squares but 
treat the
cell averages as pointwise values at the cell centers. Here again
the second derivative terms are ignored.  We refer to this as ``pointwise
quadratic reconstruction". Note that this is not a
second order accurate gradient, since the centroid value differs by $O(h^2)$
from the pointwise value at the centroid. Nevertheless, this is frequently
done since it is easier to implement. This is particularly true  for SRD
neighborhoods with irregular shapes.

{
\small
\begin{table}[h]
\centering
	\hspace*{-.3in}
	\subfloat[$L_1$ volume errors. \label{tab:ex1_L1vol}]{
		\begin{tabular}{|l|c|l|l|l|}
			\hline
			$h$ & $N_x, N_y$ & 1st order grad. & ptwise grad. & 2nd order grad.   \\
			\hline
			.5297 & 27 & 6.75e-3 & 2.76e-3 & 2.45e-3 \\
			\hline
			.2648 & 54  & 1.78e-3 (3.8)  & 6.38e-4 (4.3) & 4.71e-4 (5.2) \\
			\hline
			.1324 &108 & 3.63e-4 (4.9)  & 1.53e-4 (4.2) & 1.21e-4 (3.9) \\
			\hline
			.662e-2 & 216 & 6.52e-5 (5.6)  & 3.65e-5  (4.2) & 2.98e-5 (4.1) \\
			\hline
			.331e-2 & 432 & 1.40e-5 (4.7)  & 8.85e-6  (4.1) & 7.86e-6 (3.8) \\
			\hline
			.166e-2 & 864 & 2.68e-6 (5.2)  & 2.15e-6  (4.1) & 2.04e-6 (3.9) \\
			\hline \hline
		\end{tabular}
	}
	\quad
	\vspace*{.2in}
	
	\hspace*{-.3in}
	\subfloat[$L_1$ boundary errors. \label{tab:ex1_L1bndry}]{
		\begin{tabular}{|l|c|l|l|l|l|}
			\hline
	    $h$ &  $N_x, N_y$   & 1st order grad.  & ptwise grad. &  2nd order grad.    \\
	\hline
     .5297 & 	27  &   6.84e-02 &  3.31e-02  & 2.37e-2\\
	\hline
     .2648 &	54  &   2.83e-02 (2.4) &  1.21e-02 (2.7)  & 8.183-3  (2.9) \\
	\hline
      .1324 &	108 &   1.03e-02  (2.8)&  4.69e-03 (2.6) & 3.453-3 (2.4) \\
	\hline
     .662e-2 &	216 &   3.65e-03  (2.8) &  1.82e-03 (2.6) & 1.38e-3 (2.5) \\
			\hline
    .331e-2 &	432 &   1.24e-03  (3.0)  &  7.18e-04 (2.6) & 6.15e-4 (2.2) \\
			\hline
    .166e-2 &	864 &   3.96e-04  (3.1)  &  2.85e-04 (2.6) & 2.58e-4 (2.4) \\
			\hline \hline
		\end{tabular}
	}
	
	\caption{\sf L1 norm of the error in the domain volume and along the boundary,
        for the supersonic vortex example. The error using first order accurate
        gradients, pointwise quadratic, and fully second order accurate gradients is shown. 
        Pointwise quadratic gradients have half of the error than first order gradients. 
        Truly second order accurate gradients are even better,
        especially on coarser grids. The convergence rate is given in
        parentheses. \label{tab:ssv}}
\end{table}
}

We measure the $L_1$ norm of the error in the volume,   $\sum_{i,j} \,
V_{ij} \lvert e_{ij } \rvert$ and at the boundary, $ \sum_{{i,j} \in \partial \Omega} \, A_{ij} \lvert e_{ij
} \rvert$.  Here
$e_{i,j}$ is the error at the cell centroid, and $A_{ij}$ is the length of the 
boundary segment in cut cell $(i,j)$.
These are given in Table \ref{tab:ex1_L1vol} and \ref{tab:ex1_L1bndry}. 
For easier comparison with other papers we also give the mesh size 
to a few digits. In  a cut cell mesh there are both solid and
flow cells, so the total number of cells is larger than the number of
flow cells, and not as easy to compare as the mesh spacing.
All interior cells use the same evolution scheme and gradients, so
the difference in the errors is solely due to the irregular cell
reconstruction scheme.

\begin{figure}[h]
	\begin{center}
		\includegraphics[height=3.0in]{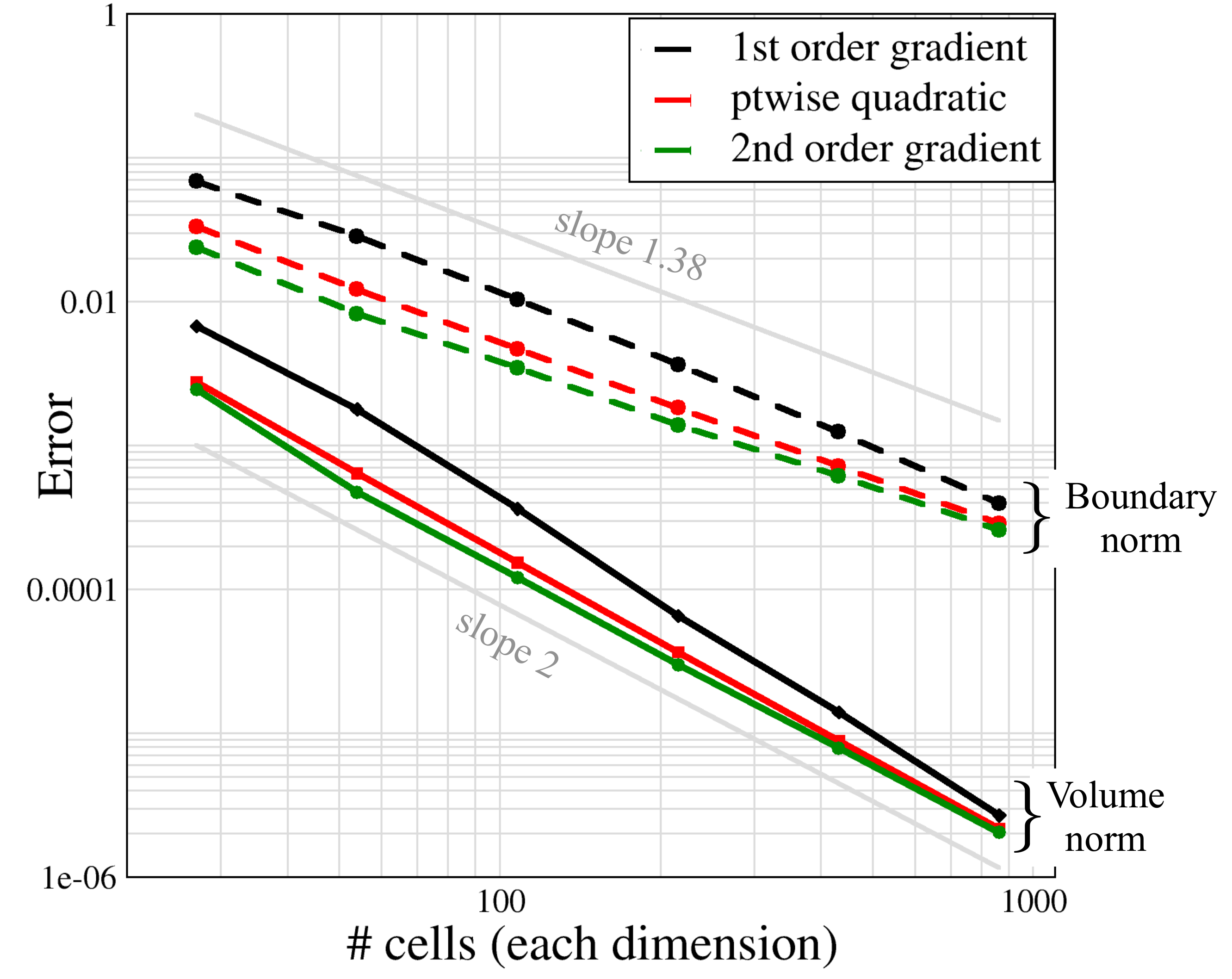}
		\caption{\sf Convergence in the L1 norm of the error in the entire 
			domain  (solid line), and along the boundary (dashed line).
			The reference line has slope 2 next to the domain error. The convergence
			rate at the boundary is between 1.38 and 1.5
			along the boundary.
			\label{fig:ssv}}
	\end{center}
\end{figure}

Note that the error at the cut cells is larger than in the volume, and has a
lower convergence rate. Since the number of cut cells grows only
linearly with refinement, the $L_1$
accuracy in the entire flow field is still second order.  At the boundary,
Richardson extrapolation shows that the convergence rate seems to be
between 1.38 and 1.5.
This has also been found in other cut cell studies
\cite{KB:2006,nemec_tm14}, and is due to the irregularity of the
difference scheme at the cut cells, and the curved boundary. 
The next example will show the same behavior using local time stepping,
which does not need SRD stabilization. It has the same convergence rate,
showing that this decrease is
not due to SRD.

Figure \ref{fig:ssv} and Table \ref{tab:ssv}  show that using pointwise 
quadratic cut cell gradients is roughly
a factor of 2 more accurate on coarser grids.  The more complicated second order
accurate gradient is even more accurate, especially on coarser grids.  Ultimately the gradient
error is reduced, and the error curves for the first and second order 
accurate gradients approach each other.

We also use this example to compare the effect of state redistribution versus 
marching to steady state using local time-stepping without SRD.  
Table  \ref{tab:ssv2} shows a comparison of the error in the converged
solution using local time stepping (LTS)  without SRD, and the error with full 
time steps and
SRD stabilization, both using  first order accurate gradients. 
The convergence rates and errors are essentially identical, showing that 
SRD does not degrade the computed solution with too much diffusion due
to the merging neighborhoods. This holds across all the other
gradient formulations too.

{
\small
\begin{table}[h]
\centering
 	\begin{tabular}{|l|c|l|l||l|l|} \hline
 		$h$ & $N_x ,N_y$ & \multicolumn{2}{|c|} {Volume Error} & \multicolumn{2}{|c|}{Boundary Error} \\ 
                \hline
 		    &            & {LTS (no SRD)} & SRD  & LTS (no SRD)  & SRD  \\ \hline
 			.5297 & 27 & 6.09e-3  &  6.75e-3   &  6.75e-02       &  6.84e-02 \\
 			\hline
 			.2648 & 54  & 1.67e-3  (3.6)  & 1.78e-4 (3.8)  &  2.79e-02  (2.4) &  2.83e-02 (2.4) \\
 			\hline
 			.1324 &108 & 3.41e-3  (4.9)  & 3.63e-4 (4.9)   &  1.00e-02  (2.8) &  1.03e-02  (2.8)\\
 			\hline
 			.662e-2 & 216 & 6.62e-5  (5.2)  & 6.52e-5 (5.6)  &  3.51e-03  (2.9) &  3.65e-03  (2.8)\\
 			\hline
 			.331e-2 & 432 & 1.34e-5  (4.9)  & 1.40e-5 (4.7)  &  1.21e-03  (2.9) &  1.24e-03  (3.0)  \\
 			\hline
 			.166e-2 & 864 & 2.59e-6  (5.2)  & 2.68e-6 (5.2)  &  3.78e-04  (3.2) &  3.96e-04  (3.1)  \\
 			\hline \hline
 	\end{tabular}
 	\caption{\sf Comparison of errors for supersonic vortex problem 
        using local time stepping, which does not use SRD, 
        and regular time stepping with SRD. 
        The errors are very similar, showing that SRD does not 
        degrade the solution with too much diffusion. \label{tab:ssv2}}
\end{table}
}

\subsection{Shock Reflection from  Cylinder}
Next we demonstrate the method using a Mach 2
shock diffracting around a circular cylinder. 

The shock will meet the cylinder and reflect at the cut cells, which are
at all angles in the mesh. This
example will demonstrate the robustness of SRD, and demonstrate the
smoothness of the density profile around the cylinder, despite the
completely irregular mesh.

A cylinder with radius $0.15$ is centered at
(0.5,0.5), and the shock is initially located at $x = 0.2$.
For this example we compare results using the MUSCL scheme and
the Method of Lines as the base schemes, both using local Lax Friedrichs for the Riemann
solver. Both methods use second order accurate gradients. The BJ 
limiter is used to limit
both the base scheme irregular cells  and neighborhood reconstruction gradients. 
The CFL for the MUSCL scheme was 0.85, and for the MOL
scheme was 1.0 using \eqref{eq:vn1}, based on the full cell volumes.

\begin{figure}[h]
\centering
\vspace*{-.25in}
\hspace*{-.2in}
\mbox{
\includegraphics[width=0.50\linewidth,trim=10 10 200 10,clip]{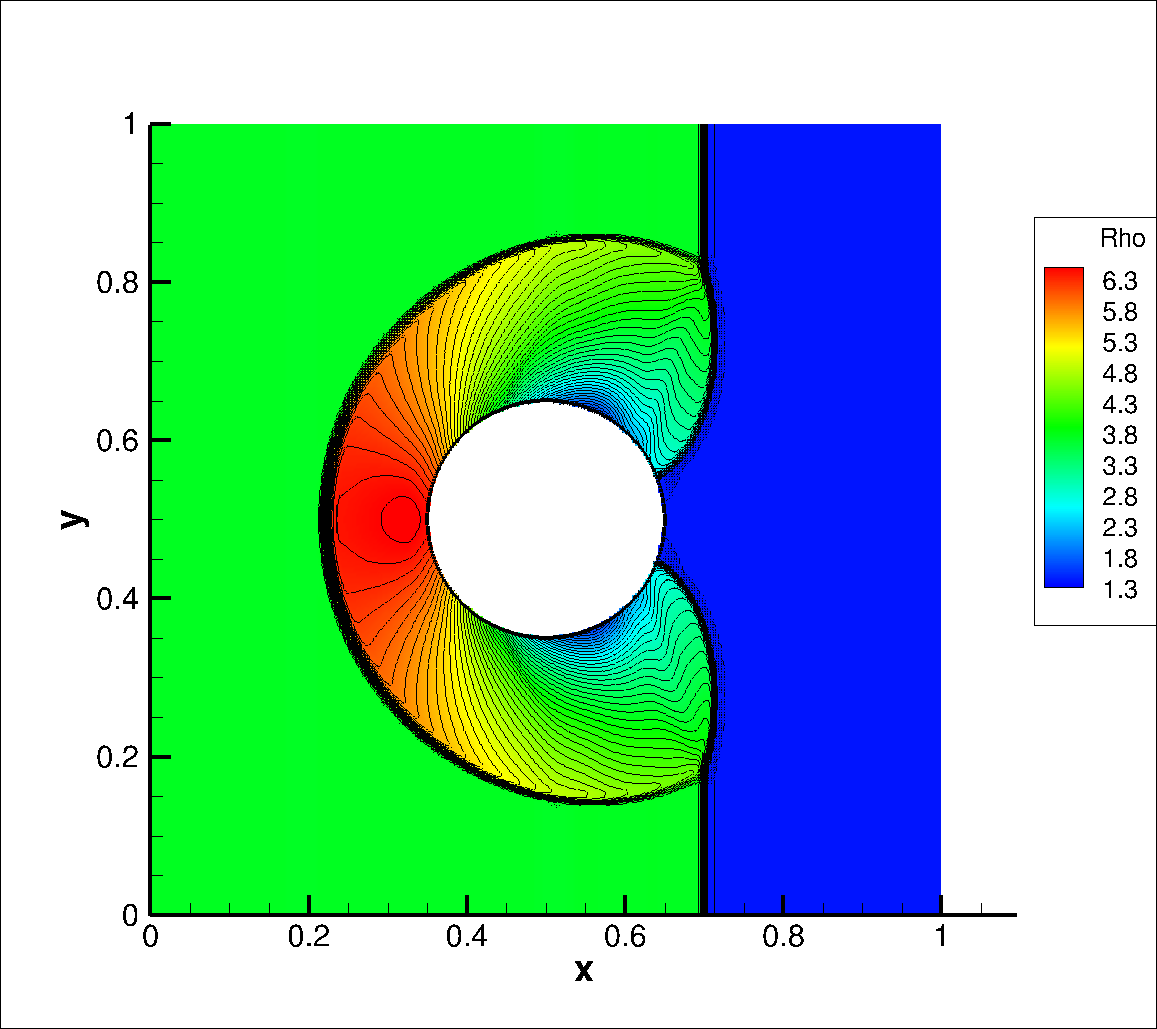}
\hspace*{.15in}
\includegraphics[width=0.50\linewidth,trim=10 10 200 10,clip]{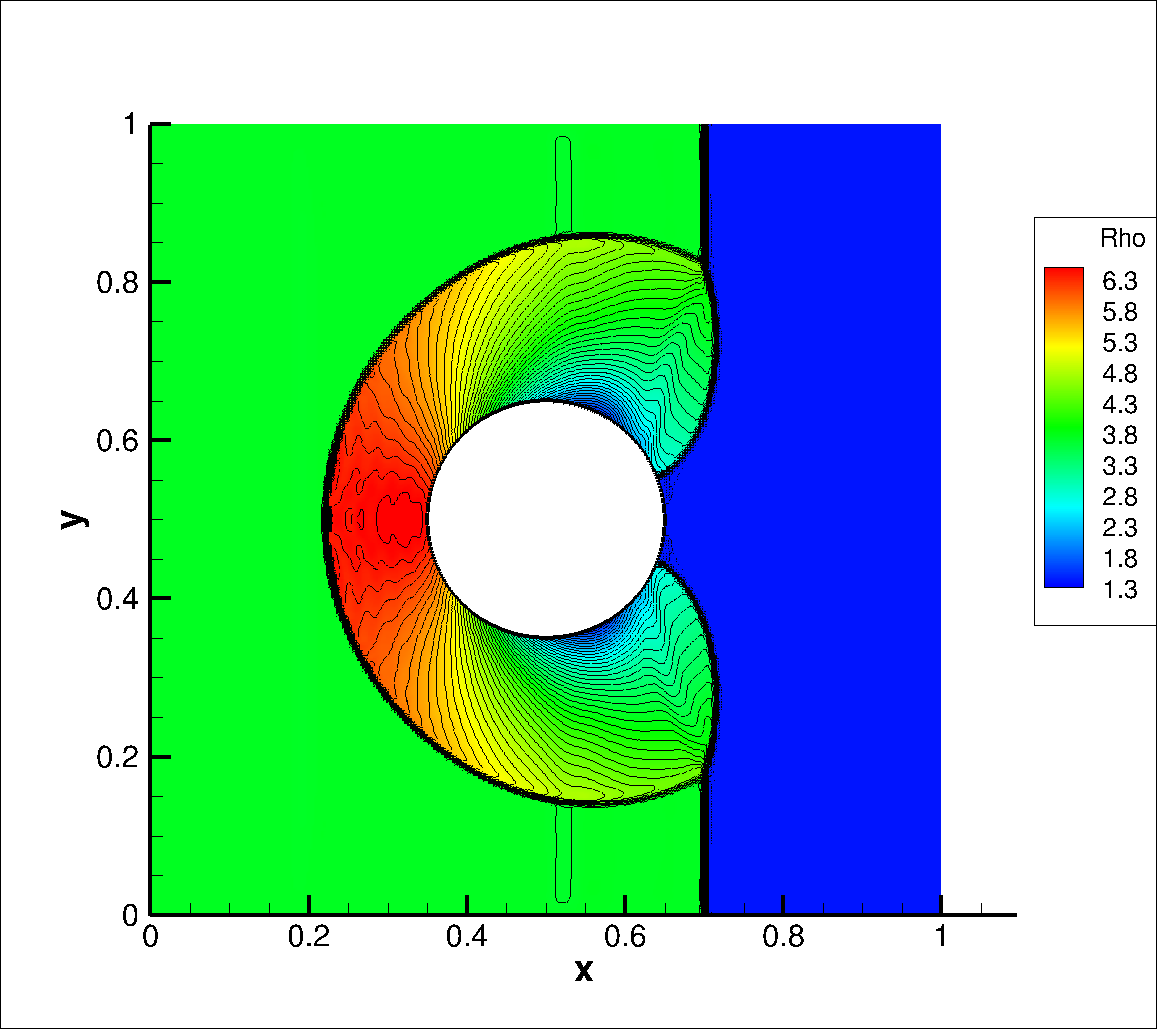}
}
\caption{\sf Density profile of Mach 2 shock reflection around a cylinder,
at time $T = 0.25$.  Left computation uses MUSCL, right uses MOL. 
There are 52 contours between 1.3 and 6.5.
Both schemes robustly compute the solution. 
The subsonic region in front of the cylinder is better behaved with MUSCL.
\label{fig:cyl1}}
\vspace*{-.1in}
\end{figure}

Figure \ref{fig:cyl1} left shows the solution density from MUSCL, and
right from MOL, at time $T=0.25$. 
Both grids use 302 cells in each
direction, and the domain is  $[0.0,1.00001,]$ by  $[0.0, 1.0]$, again to
prevent mesh degeneracies.
There are 416 cut cells
around the cylinder; 160 of the cut cells had volume fractions less than
0.5 and were stabilized with SRD.  The smallest volume fraction was 1.17e-4.   
For comparison, this is also the time shown in \cite{mjb-hel-rjl:hbox2}. 
Here and in \cite{mjb-hel-rjl:hbox2}, the solution in front of the cylinder 
where the
maximum density is located, behaves better with the one-step methods than with
MOL.  The method of lines solution is smoother around the boundary than our
MUSCL variant. 

\begin{figure}[h]
\centering
\hspace*{-.2in}
\includegraphics[height=2.7in]{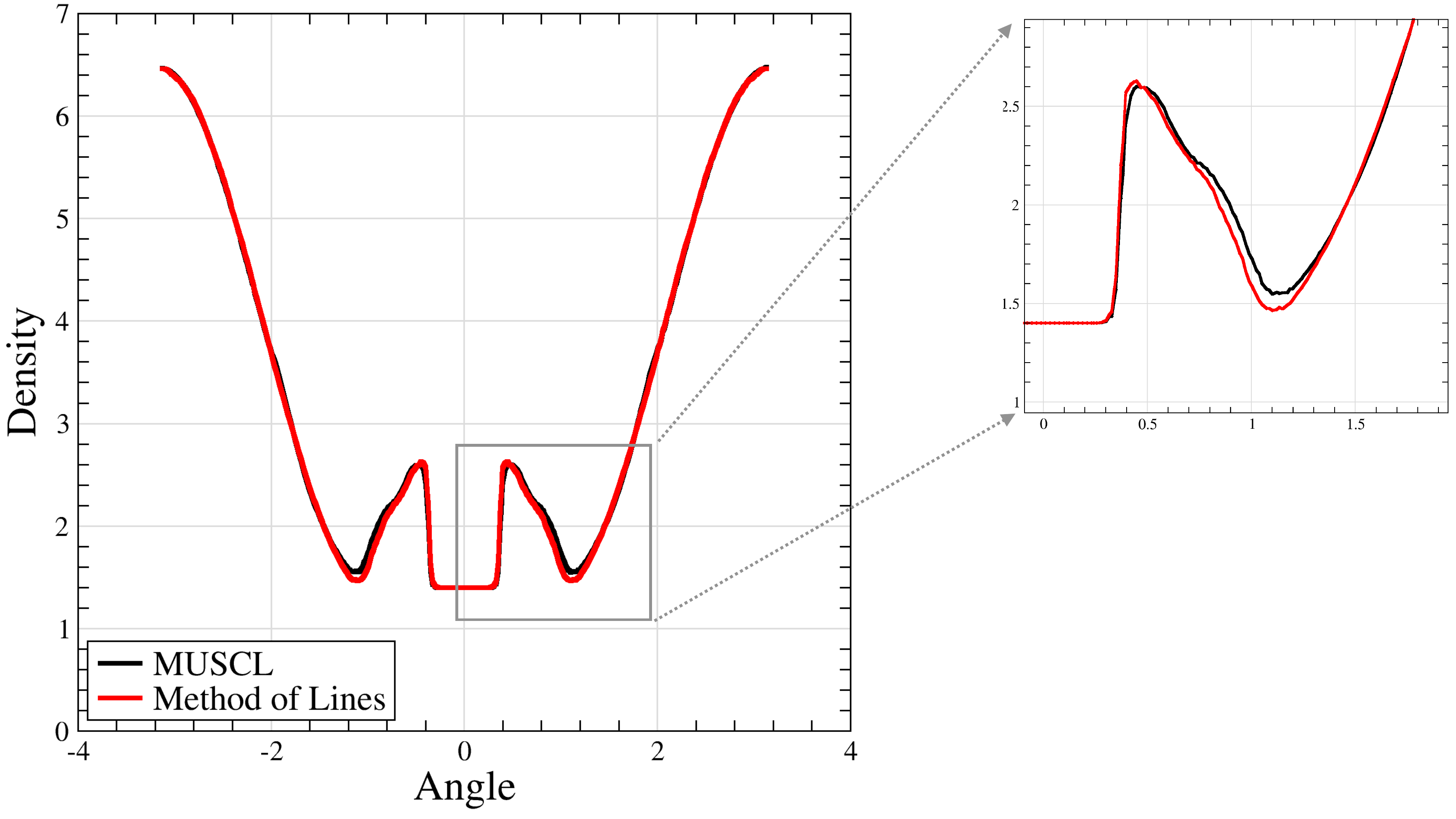}
\hspace*{.3in}
\caption{\sf The density profile  around the cylinder for the MUSCL and MOL schemes 
at time $T=0.25$.  The methods give somewhat different profiles at the
local minima, but the shock is  located at the same position.
Both methods  are remarkably smooth considering the
irregularity of the cut cell mesh.}
\label{fig:cylbndry}
\end{figure}

Figure \ref{fig:cylbndry}
shows the density profile from both schemes
taken along the cylinder.
For this plot, the cut cell variable is
reconstructed to the midpoint of the cylinder line segment in each  cell.  
The zoom shows the difference more
clearly.

\clearpage

\subsection{Double Mach Reflection problem}\label{sec:dm}
We present one final standard test case for embedded
boundary meshes. We reflect a Mach 10 shock obliquely over a wedge, 
where the shock and wall form a $60^{\circ}$ angle.  
The problem domain is $[0,3.0]\times[0,1.75]$, with an angled wall 
passing through the point $(1/6,0)$. 
We use the MUSCL scheme, with the interior slopes 
computed using monotonized central differences. The irregular cells and
SRD neighborhoods use second order accurate 
gradients, and limit using BJ.
The CFL is 0.85.

\begin{figure}[H]
\centering
\hspace*{-.2in}
\includegraphics[width=1. \linewidth]{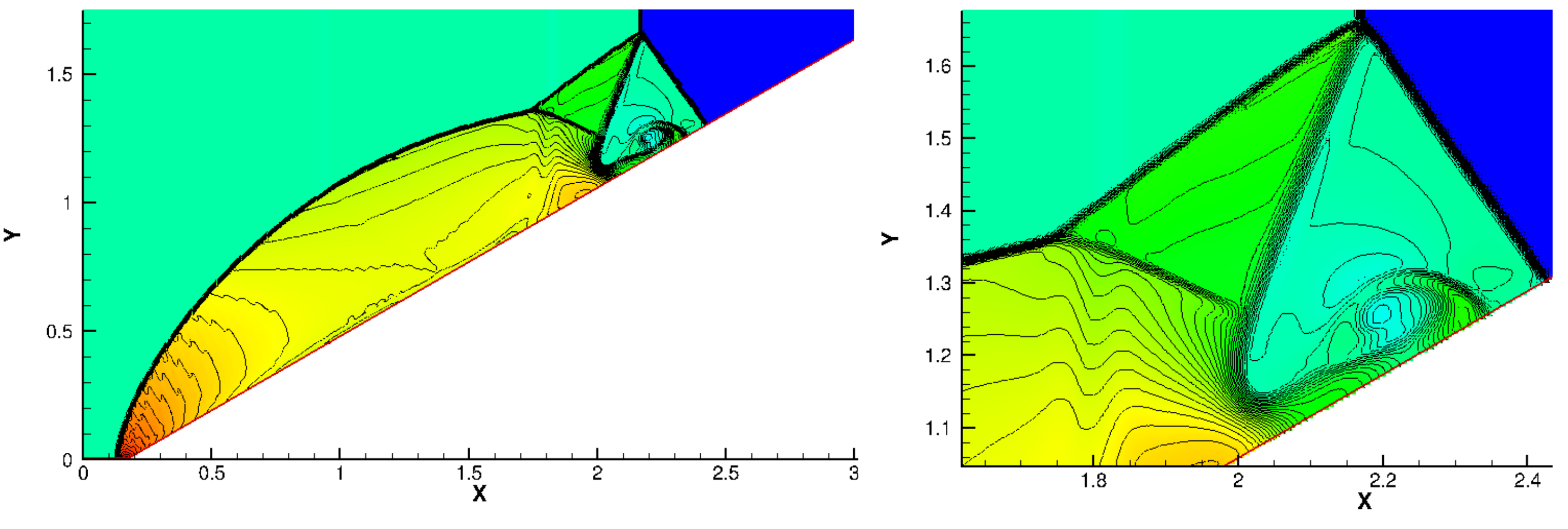}
\hspace*{-.2in}
\includegraphics[width=1. \linewidth]{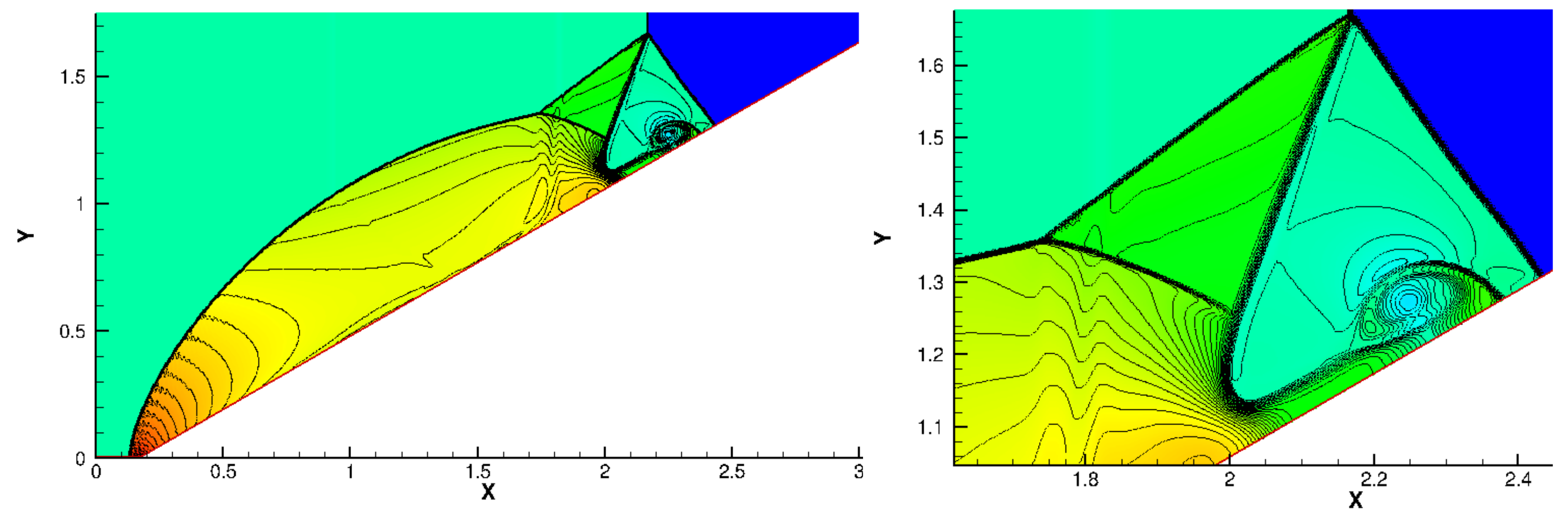}
\caption{\sf Density plot and isolines of a Mach 10 shock impinging 
obliquely on a wedge at time $T = 0.2$ where $\Delta x = \Delta y =
1/240$ (top) and $1/480$ (bottom).
The solution on the full domain is shown adjacent to a zoom of the Mach 
stem region.  Sixty isolines between 1.39 and 21.0 are drawn.}\label{fig:dm}
\end{figure}

The solution to this problem is a complex, self-similar reflection pattern 
composed of incident and reflected shocks and contact 
discontinuities \cite{WOODWARD1984115,rkdg5}.  
The contact discontinuities are an unstable feature of the solution that can be 
difficult to resolve correctly, especially in the neighborhood of the 
reflecting boundary where the carbuncle phenomenon can occur \cite{KEMM2018596}.

\begin{figure}[h]
\centering
\includegraphics[width=.7\linewidth]{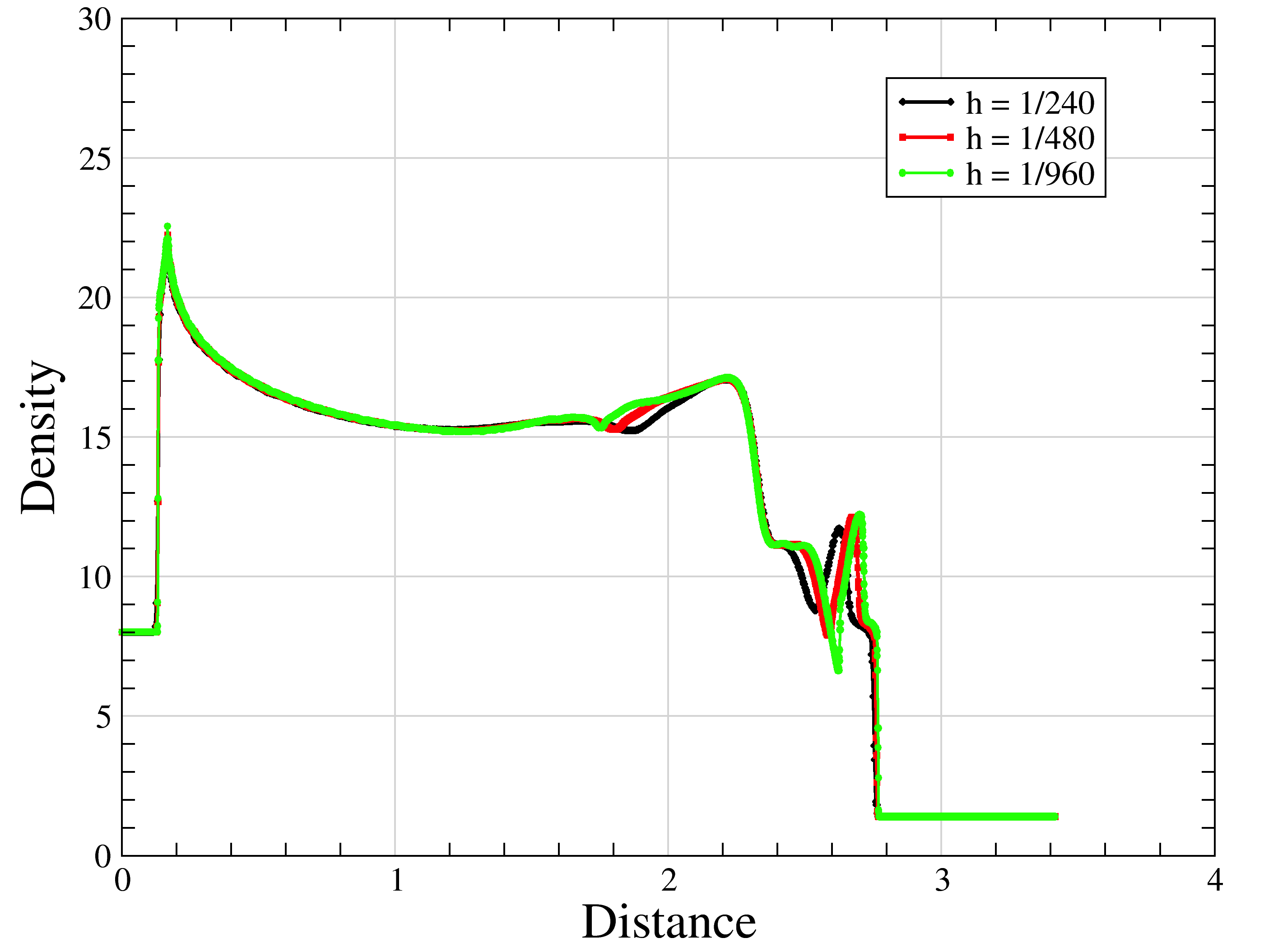}
\caption{\sf Density plotted along the wall, as a function of distance from
the lower left corner of the domain.  A finer run is included to show
the trend. The shocks are in the same place on all grids, but the 
unstable contact is slightly shifted.
} \label{fig:wedgeBndry}
\end{figure}

The solution at the final time $T = 0.2$ is plotted in Figure \ref{fig:dm}, 
where the grid resolution is $\Delta x = \Delta y = 1/240$.
We obtain qualitatively comparable results to those in \cite{rkdg5} on 
the same grid resolution. We also show for comparison the results using 
$\Delta x = \Delta y = 1/480$, also done in \cite{rkdg5}. Again, the
improvement is very similar. Finally, in Figure \ref{fig:wedgeBndry} we
show the density along the boundary for three resolutions. Despite the
irregularity of the cut cells, the solution is very smooth. The smallest
cut cell in the coarser grid has volume fraction 1.65e-6. On the finest
grid the smallest volume fraction is 2.95e-7. 

\section{Conclusions}\label{sec:conc}
We have presented a  state redistribution algorithm to solve the small cell problem on cut cell meshes.  It is conservative, allows for overlapping temporary merging neighborhoods so is easy to implement, and is linearity preserving.   
Numerical experiments show that on smooth problems, second order accuracy is maintained, and the solution is not degraded by the postprocessing. For problems with shocks, the scheme maintains robustness at the cut cells.
We have shown experiments using SRD on two different base schemes, but it should be 
applicable to any underlying numerical method with
cell-centered variables.

It should be straightforward to apply state redistribution to
three dimensional applications.
It should also be applicable to 
different sets of equations, e.g. incompressible flow. 
It seems clear that 
when used in conjunction with a higher order base scheme, 
state redistribution can be extended to higher order accuracy. 
We have already started doing this for 
3rd and 4th order accuracy. However, higher order methods bring 
in many new features, so we do not include that here.

\small
\subsection*{Acknowledgments} 
We thank Michael Aftosmis, Sandra May, and Marian Nemec 
for carefully reading the manuscript and many helpful suggestions.
This research was partially supported by the U.S. Department of
Energy under contract DE-FG02-88ER25053.

\bibliography{main}
\bibliographystyle{plain}

\end{document}